\documentclass[11pt]{article}
\usepackage{CJK}
\usepackage{mathrsfs}
\usepackage{amsmath}
\usepackage{amssymb,mathrsfs,amsthm,amsfonts}
\usepackage{graphicx}
\usepackage[pdfstartview=FitH]{hyperref}
\usepackage[title]{appendix}
\usepackage{color}

\providecommand{\U}[1]{\protect\rule{.1in}{.1in}}
\allowdisplaybreaks 

\newtheorem{thm}{Theorem}[section]
\newtheorem{proposition}{Proposition}[section]
\newtheorem{corollary}{Corollary}[section]
\newtheorem{lemma}{Lemma}[section]
\newtheorem{definition}{Definition}[section]
\newtheorem{Rem}{Remark}[section]
\newtheorem{con}{Condition}[section]
\def\mE{{\mathbb E}}
\def\eps{\varepsilon}
\def\]{{\Big]}}
\def\[{{\Big[}}
\def\dif{{\mathord{{\rm d}}}}
\def\bd{\begin{definition}}
\def\ed{\end{definition}}
\def\bp{\begin{proposition}}
\def\ep{\end{proposition}}
\def\bc{\begin{corollary}}
\def\ec{\end{corollary}}
\def\bx{\begin{Examples}}
\def\ex{\end{Examples}}
\def\cA{{\mathcal A}}
\def\cB{{\mathcal B}}
\def\cC{{\mathcal C}}

\def\cF{{\mathcal F}}

\def\cL{{\mathcal L}}

\def\cZ{{\mathcal Z}}

\def\mE{{\mathbb E}}
\def\mF{{\mathbb F}}

\def\mN{{\mathbb N}}

\def\mP{{\mathbb P}}

\def\mR{{\mathbb R}}

\def\ba{{\begin{align}}
\def\ea{\end{align}}}

\def\geq{\geqslant}
\def\leq{\leqslant}

\def\^{\widehat}

\def\ba{\begin{aligned}}
\def\ea{\end{aligned}}
\def\be{\begin{equation}}
\def\ee{\end{equation}}
\def\ben{\begin{align*}}
\def\enn{\end{align*}}

\makeatletter

\newcommand{\Rmnum}[1]{\expandafter\@slowromancap\romannumeral #1@}
\makeatother

\numberwithin{equation}{section}

\title{\bf{Well-posedness of  Stochastic 2D Hydrodynamics type Systems with Multiplicative L\'evy Noises}}

{\author{
{\ Xuhui  Peng$^\dag$ } \\
{\em\small  MOE-LCSM, School of Mathematics and Statistics, Hunan Normal University}\\
{\em\small  Changsha, Hunan 410081, P. R. China} \\
 {\em\small Email address: xhpeng@hunnu.edu.cn}\\
 {\ Juan Yang$^\ddag$ } \\
{\em\small  School of Science, Beijing University of Posts and Telecommunications}\\
{\em\small  Beijing 100876, P. R. China} \\
 {\em\small Email address: juanyang@bupt.edu.cn}\\
{\ Jianliang Zhai$^\spadesuit$} \\
{\em\small Key Laboratory of Wu Wen-Tsun Mathematics, CAS, School of Mathematical Science, }\\
{\em\small University of Science and Technology of China, Hefei 230026, P.R. China}\\
{\em\small Email address: zhaijl@ustc.edu.cn }\\
 }
 }

 \date{}

\begin{document}
\maketitle
\let\thefootnote\relax\footnotetext{$^\dag$Xuhui Peng was supported by NNSFC (No. 12071123), the Hunan Provincial Natural Science
Foundation of China (No. 2019JJ50377), the Scientific Research Project of Hunan Province Education
Department (No. 20A329) and the Construct Program of the Key Discipline in Hunan Province.}
\let\thefootnote\relax\footnotetext{$^\ddag$Corresponding Author.}
\let\thefootnote\relax\footnotetext{$^\spadesuit$Jianliang Zhai's research is supported by National Natural Science Foundation of China (NSFC) (Nos. 11971456, 11721101, 11431014,  11671372), the Fundamental Research Funds for the Central Universities
(No. WK3470000016), School Start-up Fund (USTC) KY0010000036.}
\begin{abstract}
\noindent  
This paper presents the existence and uniqueness of solutions to an abstract
nonlinear equation driven by a multiplicative noise of L\'evy type, which covers many hydrodynamical models including 2D Navier-Stokes equations,
2D MHD equations, the 2D Magnetic Bernard problem, and several Shell models of
turbulence. In the existing literature on this topic, besides the classical Lipschitz and linear growth conditions, other assumptions, which might be untypical,  are also required on the coefficients of the stochastic perturbations.
The goal of this paper is to get rid of these untypical assumptions.
{Our assumption on the coefficients of stochastic perturbations is new even for the Wiener cases, and in some sense, is shown to be quite sharp.}
A new cutting-off argument and energy estimation procedure play an important role in establishing the existence and uniqueness under this assumption.


\vspace{1mm}
\end{abstract}
\section{Introduction}
\label{}


Stochastic partial differential equations (SPDEs) driven by jump-type noises such as L\'evy-type or Poisson-type perturbations are drastically
different from those driven by Wiener noises, due to the presence of jumps, concerning the well-posedness, the Burkholder-Davis-Gundy inequality, the Girsanov theorem, the time regularity,
the ergodicity, irreducibility, mixing property and other long-time behaviour of the solutions. In general, all the results and/or techniques available for the SPDEs with Gaussian noise are not always suitable for the treatment of SPDEs with L\'evy
noise and therefore we require new and different techniques. For more details, please refer to \cite{BICHTELER, BGIPPZ, FHR, Jacod-Shiryaev, KPS, Liu Zhai 2012, Liu Zhai 2016, PZ} and the references therein.

In this paper, we are concerned about the well-posedness of SPDEs with multiplicative L\'evy noise and
the reader is referred to
\cite{Pe15} for a thorough introduction of SPDEs with L\'evy noises.
There are extensive results on the well-posedness of SPDEs with Gaussian noise.
 Here we only list several of them (see \cite{Bouard Debussche01, Bouard Debussche02, Brzezniak Millet, Chueshov Millet, Liu R, LR, Liu Zhao} and the references therein). The standard assumptions
on the coefficients of the Wiener noises are the classical Lipschitz and one sided linear growth assumptions. The main approaches to solving  SPDEs with Gaussian noises are local monotonicity
arguments combined with Galerkin approximation methods, the cutting-off and the Banach fixed point theorem (see e.g., \cite{Liu R, LR} and \cite{Brzezniak Millet, Bouard Debussche01} respectively). But using the same idea to the L\'evy case, one needs to assume
other conditions on the coefficient of L\'evy noise (see e.g., \cite{liuwei, BHZ}). In fact, besides
the classical Lipschitz and one sided linear growth assumptions, the solvability of
SPDEs with L\'evy noise always requires other assumptions on the coefficients of the stochastic perturbations in the existing literature.
And we will give more details below.

To solve SPDEs with L\'evy noises, one natural approach is based on approximating the Poisson random measure $\eta$ by a sequence of Poisson random measures $\{\eta_n\}_{n\in\mathbb{N}}$ whose
intensity measures are finite. Dong and Xie \cite{DX2009} used this approach to establish the well-posedness of the strong solutions in probabilistic  sense for 2D stochastic Navier-Stokes equations with L\'evy noise. However, besides (H1) and (H2) with $L_2=L_5=0$ in Condition \ref{con G} below, this method needs the following assumption on the control of the ``small jump": for some $k>0$,
\begin{eqnarray*}
  \sup_{|v|\leq k}\int_{\|z\|_{\cZ}\leq \delta} |G(t,v,z)|^2\nu(\dif z)\rightarrow 0, \text{  as } \delta\rightarrow 0.
\end{eqnarray*}

Another approach is the local monotonicity method combined with the Galerkin approximation. This approach is suitable for treating SPDEs with Wiener
noise (see e.g. \cite{Liu R}). However, applying this approach to the L\'evy cases (see \cite{liuwei, BHZ}), the following assumption is essential: for some $p>2$,
\begin{eqnarray*}
  \int_{\cZ} |G(t,v,z)|^p \nu(\dif z)\leq K(1+|v|^p).
\end{eqnarray*}
In particular, applying their framework to 2D stochastic Navier-Stokes equations, they need to assume that
\begin{itemize}
 \item[(J1)](Lipschitz) $$
             \|\Psi(t,v_1)-\Psi(t,v_2)\|_{\cL_2}^2
             +
             \int_\mathcal{Z}|G(t,v_1,z)-G(t,v_2,z)|^2\nu(dz)\leq L_1|v_1-v_2|^2+L_2\|v_1-v_2\|^2;
            $$

\item[(J2)](Growth) $$
             \|\Psi(t,v)\|^2_{\cL_2}
             +
             \int_{\mathcal{Z}}|G(t,v,z)|^2\nu(dz)\leq L_3+L_4|v|^2+L_5\|v\|^2;
            $$
\item[(J3)] $L_2\in[0,2)$ and $L_5\in[0,\frac{1}{2})$;

\item [(J4)] $\int_{\cZ} |G(t,v,z)|^4 \nu(\dif z)\leq K(1+|v|^4)$.
\end{itemize}
See (H1)-(H4), Theorem 1.2, and Page 292 in \cite{liuwei}. 
They can not deal with the case of $\Psi(t,v)\equiv 0$
and $G(t,v,z)=\theta z\nabla v$ with $\theta\not=0$ (see (\ref{p-1})). To the best of our knowledge, for the Wiener case, i.e., $G(t,v,z)\equiv0$, these assumptions (J1)-(J3) are the best in the existing literature.
We also refer to \cite{CWW19, DZ11, Mo13, Mo14, SS12} and the references therein, in which the existence of martingale solutions for SPDEs with L\'evy noises were established, and in these papers, the Galerkin approximation method had been applied, and some untypical
assumption like (J4) was needed.
Consequently, this motivates a unified approach to get rid of untypical assumptions
on the coefficients of the stochastic perturbations, and also makes it possible to
cover a wide class of mathematical coupled models from fluid dynamics. Our unified approach is based
on an abstract stochastic evolution equation (\ref{p-1}). It covers many hydrodynamical models (see Remark \ref{Rem 1}).
Our assumption on the coefficients of stochastic perturbations (see Condition \ref{con G}) is new even for the Wiener cases, and in some sense, is shown to be quite sharp (see Remark \ref{Rem 2} for more details).
We apply localization arguments and fixed point methods. A new cutting-off argument and energy estimation procedure play an important role.


{Together with Zdzis{\l}aw Brze\'{z}niak}, Peng and Zhai introduced a new cutting-off argument in \cite{BPZ}, showing that there exists a unique  strong solution in probability sense to stochastic 2D Navier-Stokes equations with L\'evy noises under (H1) and (H2) with $L_2=L_5=0$ in Condition \ref{con G}.
However, this method can not be suitable to deal with the case of $L_2, L_5\in(0,2)$. In this paper, we employ
a slightly different localization methods, and using a similar cutting-off argument, we establish  
{finer a  priori} estimate of $I_n$ than that in \cite{BPZ} (see pages 9--15 in \cite{BPZ} and Lemma \ref{lemma 2} in this paper).
And a new energy estimation procedure has been introduced to obtain the Lipchiz property of $\{y_n, n\in\mathbb{N}\}$ (see Propositions \ref{Th 01} and \ref{thm 02}). The whole programme is technical and highly non-trivial. We refer to Remark 2.2 (Page 390) in \cite{Chueshov Millet} for the reason why we consider the case $L_2, L_5\in(0,2)$.

\vskip 0.2cm
We should mention that the existence and uniqueness of strong solution in PDE sense to several stochastic hydrodynamical systems with L\'evy noises was established in \cite{BHR}. They also used some localization arguments and fixed point methods, which are different from this paper. Their approach  required that, for any $x,y\in V$ and $q=1,2$, there exists a constant $\ell_q>0 $ such that
\begin{eqnarray*}
  \int_{\cZ}\|G(t,x,z)-G(t,y,z)\|^{2q}\nu(dz)\leq \ell_q^q \|x-y\|^{2q}.
\end{eqnarray*}
In \cite{BPZ},  
a similar idea as in this paper has been used to prove the existence and uniqueness of strong solution in PDE sense to stochastic 2D Navier-Stokes equations with L\'evy noises without the above assumption with $q=2$. And we believe that our method in this paper and \cite{BPZ} can be used to deal with other SPDEs and PDEs.

\vskip 0.2cm
The layout of the present paper is as follows. In Section \ref{Sec 2}, we introduce the abstract stochastic evolution
equation that our result is based on, and our main result. Section \ref{pp-1} is devoted to the proof of our main result.

\vskip 0.2cm

\section{Preliminaries and Main Result} \label{Sec 2}
Suppose that $(\Omega,\cF,\mF=\{\cF_t\}_{t\geq0},\mP)$ is a filtered probability space satisfying the
  usual conditions.
 Let $(\cZ,\cB(\cZ))$ be  a measurable space, $\nu$ be a given $\sigma$-finite measure on $\cZ$,
 $ \nu(\{0\}) = 0$. Let $\eta:\Omega\times\cB(\mR_+)\times \cB(\cZ)\rightarrow\bar{\mN}$
 be a time homogeneous Poisson random
measure on $(\cZ,\cB(\cZ))$ with intensity measure $\nu$. We write
$$\tilde{\eta}([0,t]\times O)=\eta([0,t]\times O) -t\nu(O),\ \ t\geq0,\ O\in\mathcal{B}(\mathcal{Z})$$
for the compensated Poisson random measure associated with $\eta$.
Let
$\{ W(t)\}_{t\geq 0}$
be a  $K$-valued cylindrical Wiener process
on $(\Omega,\cF,\mF,\mP)$, where $K$ is a separable Hilbert space.

The aim of this paper is to study the well-posedness of the abstract evolution equation given by
\begin{eqnarray}\label{p-1}
\left\{
  \begin{split}
    & \dif u(t)+\cA u(t)\dif t+B(u(t),u(t))\dif t
     \\ & = f(t)dt+ \int_{\mathcal{Z}}G(t,u(t-),z)\widetilde{\eta}(dz,dt)+\Psi(t,u(t))d W(t),
    \\
    & u(0)=u_0\in H,
  \end{split}
  \right.
\end{eqnarray}
where $H$ is a separable Hilbert space, and $\cA$ is a (possibly unbounded) self-adjoint positive linear operator on $H$.
We shall denote the scalar product and the norm of $H$ by $\langle\cdot,\cdot \rangle$ and $|\cdot|$ respectively.
Set $V=\rm{dom} (\cA^{1/2})$  equipped with norm $\|x\|:=|\cA^{\frac{1}{2}}x|,x\in V.$
The operator $B:V\times V\rightarrow V'$ is a continuous mapping, where $V'$ is the dual of $V.$
With a slightly abuse of notation, the duality between $V'$ and $V$ is also denoted by $\langle f,v\rangle$ for $f\in V'$ and $v\in V$, whose meaning should be clear from the context. The coefficients of the stochastic perturbations  $G$ and $\Psi$ are
measurable functions, satisfying certain conditions specified later.

\vskip 0.2cm

Our condition on the operator $B$ is the following.
\begin{con}\label{p-2}
Assume that $B:V\times V\rightarrow V'$   is a continuous bilinear mapping satisfying the following conditions:
\begin{itemize}
\item[(B1)](Skewsymmetricity of $B$)
\begin{eqnarray}\label{p-3}
  \langle B(u_1,u_2),u_3\rangle=-\langle B(u_1,u_3),u_2\rangle, \text{ for all } u_1,u_2,u_3\in V;
\end{eqnarray}

\item[(B2)] there exists a reflexive and separable Banach space $(Q, |\cdot|_Q)$ and a constant $a_0>0$ such that
    \begin{eqnarray}
      && V\subset Q \subset H,
      \\ \label{p-4} && |v|_Q^2\leq a_0|v|\cdot \|v\|, \text{  for all } v\in V;
    \end{eqnarray}
\item[(B3)] there exists a constant $C>0$ such that
\begin{eqnarray}
\label{p-5}
  |\langle B(u,v), w\rangle|\leq C|u|_Q \|v\| |w|_Q,
  \text{  for all } u,v,w\in V.
\end{eqnarray}
\end{itemize}
\end{con}

\begin{Rem}
(\ref{p-3}) implies that \begin{eqnarray}\label{eq PB}
  \langle B(u_1,u_2),u_2\rangle=0, \text{ for all } u_1,u_2\in V.
\end{eqnarray}
\end{Rem}

\begin{Rem}\label{Rem 1}
This type of  abstract evolution  equation (\ref{p-1}) covers stochastic 2D Navier-Stokes equations,
2D stochastic Magneto-Hydrodynamic equations, 2D stochastic Boussinesq model
for the B\'{e}nard Convection, 2D stochastic Magnetic B\'{e}rnard problem, 3D stochastic
Leray $\alpha$-Model for Navier-Stokes equations and several stochastic Shell models of turbulence. For more details,  we refer to \cite{BHZ},  \cite[Section 2.1]{CM}.
\end{Rem}

Now, we give the definition of the solution to (\ref{p-1}).

\begin{definition}
An $H$-valued $\rm c\grave{a}dl\grave{a}g$ $\mF$-adapted process $\{u(t)\}_{t\in [0,T]}$ is called a solution of (\ref{p-1}) if the following conditions are satisfied,

 (S1)  $u\in D([0,T],H)\cap L^2([0,T],V)$, $\mathbb{P}$-a.s.,
where $D([0,T],H)$ denotes all of the $\rm c\grave{a}dl\grave{a}g$ functions from $[0,T]$ into $H$ equipped with the Skorohod topology.

  (S2) the following  equality holds for every $t\in [0,T]$, as an element of $V'$, $\mP\text{-}$a.s.
        \begin{eqnarray*}\label{eq S2}
    u(t)&=&u_0-\int_0^t\mathcal{A}u(s)ds - \int_0^tB(u(s),u(s))ds+\int_0^t f(s)ds\nonumber
    \\            &&+\int_0^t\int_{\mathcal{Z}}G(s,u(s-),z)\widetilde{\eta}(dz,ds)+\int_0^t\Psi(s,u(s))d W(s).
           \end{eqnarray*}
\end{definition}
An alternative version of condition (S2) is to require that for every $t\in [0,T]$, $\mP\text{-}$a.s.
\begin{eqnarray*}
   \langle u(t),\phi\rangle &=&\langle u_0,\phi\rangle -\int_0^t\langle\mathcal{A}u(s), \phi\rangle ds - \int_0^t\langle B(u(s),u(s)),\phi\rangle ds+\int_0^t \langle f(s),\phi\rangle ds
   \\  &&+\int_0^t\int_{\mathcal{Z}} \langle G(s,u(s-),z), \phi \rangle\widetilde{\eta}(dz,ds)+\int_0^t\langle \Psi(s,u(s))d W(s),
   \phi \rangle, \ \forall \phi\in V.
\end{eqnarray*}

Before presenting the main result of this paper, let us first formulate the main assumptions on the coefficients $G$ and $\Psi$.
Let us denote by $(\cL_2(K,H),\|\cdot\|_{\cL_2})$ the Hilbert space of all Hilbert-Schmidt operators from  $K$ into $H$.
\begin{con}\label{con G}
$G:[0,T]\times H\times \mathcal{Z}\rightarrow H$  and $\Psi:[0,T]\times H\rightarrow\cL_2(K,H)$  are measurable mappings.  There exist constants  $ L_i>0,\ i=1,\cdots,5$ such that, for all $t\in[0,T],\ v, v_1, v_2\in V$,
\begin{itemize}
\item[(H1)](Lipschitz) $$
             \|\Psi(t,v_1)-\Psi(t,v_2)\|_{\cL_2}^2
             +
             \int_\mathcal{Z}|G(t,v_1,z)-G(t,v_2,z)|^2\nu(dz)\leq L_1|v_1-v_2|^2+L_2\|v_1-v_2\|^2;
            $$

\item[(H2)](Growth) $$
             \|\Psi(t,v)\|^2_{\cL_2}
             +
             \int_{\mathcal{Z}}|G(t,v,z)|^2\nu(dz)\leq L_3+L_4|v|^2+L_5\|v\|^2;
            $$
\item[(H3)] $L_2,L_5\in[0,2)$.
\end{itemize}
\end{con}

Now we state our main result, whose proof is provided in section \ref{pp-1}.
\begin{thm}\label{thm solution}
Assume that \textbf{Conditions \ref{p-1} and \ref{con G}} hold.
Then for any $\cF_0$-measurable $H$-valued initial data $u_0$ satisfying $\mE [|u_0|^2]<\infty$ and $f\in L^2([0,T],V^{'})$,
there exists a unique solution  $\{u(t)\}_{t\in[0, T]}$ to problem (\ref{p-1}). Moreover,
  \begin{eqnarray*}
    \sup_{t\in [0,T]}\mE[|u(t)|^2]+\mE[\int_0^T \|u(t)\|^2\dif t]<\infty.
  \end{eqnarray*}
\end{thm}
\vskip 0.2cm

\begin{Rem}\label{Rem 2}
Theorem \ref{thm solution} is the best in the following sense. Applying $\rm It\hat{o}$'s Formula to $|u(t)|^2$ to yield
\begin{eqnarray*}
&&|u(t)|^2+2\int_0^t\|u(s)\|^2ds\\
&=&
|u_0|^2
+
2\int_0^t\langle f(s),u(s)\rangle ds
+
2\int_0^t\int_{\mathcal{Z}}\langle G(s,u(s-),z),u(s-)\rangle\widetilde{\eta}(dz,ds)\\
&&+
2\int_0^t\langle\Psi(s,u(s)),u(s)\rangle d W(s)
+
\int_0^t\int_{\mathcal{Z}}|G(s,u(s-),z)|^2\eta(dz,ds)+\int_0^t\|\Psi(s,u(s))\|^2_{\mathcal{L}_2}d s.
\end{eqnarray*}
It is reasonable to suppose that $2\int_0^t\int_{\mathcal{Z}}\langle G(s,u(s-),z),u(s-)\rangle\widetilde{\eta}(dz,ds)$ and $2\int_0^t\langle\Psi(s,u(s)),u(s)\rangle$ $d W(s)$ are local martingales, and then using a suitable stopping time technique,
we can obtain, for any $\eps>0$,
\begin{eqnarray*}
&&\mathbb{E}[|u(t)|^2]+2\mathbb{E}[\int_0^t\|u(s)\|^2ds]\\
&\leq&
\mathbb{E}[|u_0|^2]
+
2\mathbb{E}[\int_0^t\langle f(s),u(s)\rangle ds]
+
\mathbb{E}[\int_0^t\int_{\mathcal{Z}}|G(s,u(s-),z)|^2\eta(dz,ds)]\\
&&
+
\mathbb{E}[\int_0^t\|\Psi(s,u(s))\|^2_{\mathcal{L}_2}d s]\\
&\leq&
\mathbb{E}[|u_0|^2]
+
\eps \mathbb{E}[\int_0^t\|u(s)\|^2ds]
+
\frac{1}{\eps}\int_0^t\|f(s)\|^2_{V'}ds\\
&&+
\mathbb{E}[\int_0^t\int_{\mathcal{Z}}|G(s,u(s),z)|^2\nu(dz)ds]
+
\mathbb{E}[\int_0^t\|\Psi(s,u(s))\|^2_{\mathcal{L}_2}d s].
\end{eqnarray*}
The above inequality and (H2) in Condition \ref{con G} imply that
\begin{eqnarray*}
&&\mathbb{E}[|u(t)|^2]+(2-L_5-\eps)\mathbb{E}[\int_0^t\|u(s)\|^2ds]\\
&\leq&
\mathbb{E}[|u_0|^2]
+
\frac{1}{\eps}\int_0^t\|f(s)\|^2_{V'}ds
+
L_3t
+
L_4\int_0^t \mathbb{E}[|u(s)|^2  ]d s.
\end{eqnarray*}
It is easy to see that the best assumption on $L_5$ is $L_5<2$. Using a similar argument in proving the uniqueness, the best assumption on $L_2$ is $L_2<2$.
Therefore, the well-posedness  to problem (\ref{p-1}) seems to be false if 
{(H3)} in Condition \ref{con G} does not hold.

\end{Rem}


%

\section{Proof of Theorem  \ref{thm solution}}
\label{pp-1}

We start by introducing some notations and   main ideas  used  in this paper. Then we will give the proof of Theorem \ref{thm solution}.

{\bf Notations:}

In the following, we will denote $D(I;M)$ the space of all $\rm c\grave{a}dl\grave{a}g$ paths from a time interval
$I$ into a metric space $M$.

Set
$$
\Upsilon_t= D([0,t],H)\cap L^2([0,t],V).
$$
For any  $t>0$ and $y\in \Upsilon_t$ define
\begin{eqnarray*}
|y|_{\xi_t}^2&=&\int_0^t\|y(s)\|^2ds.
\end{eqnarray*}
{Let $\Lambda_t$ be the space of all $\Upsilon_t$-valued $\{\mathcal{F}_s\}_{s\in[0,t]}$-adapted processes $y$ satisfying
$$
|y|_{\Lambda_t}^2:=\sup_{s\in[0,t]}\mathbb{E} [|y(s)|^2]  +\mathbb{E} [\int_0^t\|y(s)\|^2 ds ]<\infty.
$$
}

For any $m\in \mN,$ fix a function $\phi_m:[0,\infty)\rightarrow [0,1]$ satisfying
\begin{equation*}
\left\{\begin{array}{l}
\displaystyle
\phi_m \in C^2[0,\infty),\\
L_\phi:=\sup_{ t\in[0,\infty)}|\phi_m'(t)|<\infty,\\
\phi_m(t)=1,\ \  \  \ \ {t\in[0,m]},\\
\phi_m(t)=0, \ \ \ \ \ {t\geq m+1}.
 \end{array}\right.
\end{equation*}
Here we mention that $L_\phi$ is independent of $m$.

For any $\delta>0,$ fix a function $g_\delta:[0,\infty)\rightarrow [0,1]$ satisfying
\begin{equation*}
\left\{\begin{array}{l}
\displaystyle
g_\delta \in C^2[0,\infty),\\
\sup_{t\in[0,\infty)}|g_\delta'(t)|\leq \frac{K}{\delta},\\
g_\delta(t)=1,\ \  \  \ \ {t\in[0,\delta]},\\
g_\delta(t)=0, \ \ \ \ \ {t\geq 2\delta}.
 \end{array}\right.
\end{equation*}
Here $K$ is independent on $\delta$.

{\bf Main ideas:}

{Now we introduce the main idea in this paper, which will be divided into four steps:

{\it Step 1: Cutting-off argument.} For any $m\in \mN$, $\delta>0$ and  $y_0\in \Lambda_T$, we prove that there exists a unique solution to (\ref{eq BJ00}), which is stated in Lemma \ref{Lemma 1}.

{\it Step 2: Energy estimation.} Set $y_0(t)\equiv0$, by Lemma \ref{Lemma 1}, we define  $y_{n+1}$ satisfying (\ref{XZ-1}) recursively. Thanks to $y_0(t)\equiv0$,  we can prove that there exist  $\delta_0>0$ and $T_0>0$ such that
\begin{eqnarray*}
\sum_{n=2}^\infty\Big(\mathbb{E}[\int_0^{T_0}\|y_{n+1}(s)-y_n(s)\|^2ds]
\Big)^{1/2}
+
\sum_{n=2}^\infty\mathbb{E}[\sup_{t\in[0,T_0]}|y_{n+1}(t)-y_n(t)|
]<\infty.
\end{eqnarray*}
Here $y_{n+1}$ is the solution of (\ref{XZ-1}) with $\delta$ replaced by $\delta_0$. See Propositions \ref{Th 01}
and \ref{thm 02}. In order to do this, we need  a priori estimates. See Lemmatas \ref{lemma 2} and \ref{lemma 3}.

{\it Step 3: Local existence.} By Propositions \ref{Th 01} and \ref{thm 02}, we can prove that,
for any $T>0$ and $m>0$, there exists a solution to (\ref{p-27}) on $[0,T]$. See Proposition \ref{thm 03}.  This implies the local existence of (\ref{p-1}).

{\it Step 4: Global existence.} Finally, we prove the global existence, and for the uniqueness we refer to \cite{liuwei} or \cite{BHZ}.

%
%
%
%

\vskip 0.2cm
 \textbf{ Proof of Theorem  \ref{thm solution}:}

Now we are in a position to give the details.


\begin{lemma}\label{Lemma 1} Under the same assumptions as in Theorem  \ref{thm solution},
for any $m\in \mN$, $\delta>0$ and  $y_0\in \Lambda_T$, we have
\begin{itemize}
  \item[{\rm (Claim 1)}] for any $H$-valued progressively measurable process $h=\{h(t),\ t\in[0,T]\}$ satisfying
$$
\sup_{s\in[0,T]}\mathbb{E}[|h(s)|^2]+\mathbb{E}[
\int_0^T\|h(s)\|^2ds]<\infty,
$$
there exists a unique element $\Phi^h\in \Lambda_T$ satisfying,
for every $t\in [0,T]$, as an element of $V'$, $\mP\text{-}$a.s.
\begin{eqnarray}
 \Phi^{h}(t)&=&u_0-\int_0^t\mathcal{A}\Phi^{h}(s)ds-\int_0^tB(y_0(s), \Phi^{h}(s))\phi_m(|y_0(s)|)g_\delta(|y_0|_{\xi_s})ds\nonumber
\\ && + \int_0^tf(s)ds + \int_0^{t}\int_{\mathcal{Z}}G(s,h(s),z)\widetilde{\eta}(dz,ds)+\int_0^t\Psi(s,h(s))d W(s).\nonumber
\end{eqnarray}

  \item[{\rm (Claim 2)}] there exists a unique element 
  {$y_1=\Theta^{y_0}\in \Lambda_T$} satisfying,
for every $t\in [0,T]$, as an element of $V'$, $\mP\text{-}$a.s.
\begin{eqnarray}\label{eq BJ00}
 y_1(t)&=&u_0-\int_0^t\mathcal{A}y_1(s)ds-\int_0^tB(y_0(s), y_1(s))\phi_m(|y_0(s)|)g_\delta(|y_0|_{\xi_s})ds
\\ && + \int_0^tf(s)ds + \int_0^{t}\int_{\mathcal{Z}}G(s,y_1(s-),z)\widetilde{\eta}(dz,ds)+\int_0^t\Psi(s,y_1(s))d W(s).\nonumber
\end{eqnarray}
Moreover,
\begin{eqnarray}\label{eq BJ 03}
 \sup_{t\in[0,T]}\mathbb{E}[|y_1(t)|^2]
 +\mathbb{E}[\int_0^{T}\|y_1(s)\|^2ds]
 \leq
    C\Big(\mathbb{E}[|u_0|^2],
   \int_0^T\|f(s)\|^2_{V'}ds, T\Big).
\end{eqnarray}
Here $C\Big(\mathbb{E}|u_0|^2,
   \int_0^T\|f(s)\|^2_{V'}ds, T\Big)$ is independent of $m,\delta$ and  $y_0$.
\end{itemize}

\end{lemma}

 \noindent\emph{Proof of Lemma \ref{Lemma 1}}.}\quad
We first verify (Claim 1).

For any $H$-valued progressively measurable process $h=\{h(t),\ t\in[0,T]\}$ satisfying
$$
\sup_{s\in[0,T]}\mathbb{E}[|h(s)|^2]
+\mathbb{E}[\int_0^T\|h(s)\|^2ds]<\infty,
$$
by Condition \ref{con G}, we have
\begin{eqnarray}\label{eq Z01}
&&\mathbb{E}[\int_0^T\int_{\mathcal{Z}}|G(s,h(s),z)|^2\nu(dz)ds]
+\mathbb{E}[\int_0^T\|\Psi(s,h(s))\|^2_{\mathcal{L}_2}ds]\nonumber\\
&\leq&
L_3T+L_4T\sup_{s\in[0,T]}\mathbb{E}[|h(s)|^2]
+L_5\mathbb{E}[\int_0^T\|h(s)\|^2ds]<\infty.
\end{eqnarray}

From the classical Gelerkin approximation arguments, it is easy to prove that there exists a unique $Z^h\in \Upsilon_T,\ \mathbb{P}$-a.s. such that
for every $t\in [0,T]$, as an element of $V'$, $\mP\text{-}$a.s.
$$
Z^h(t)=-\int_0^t\mathcal{A}Z^h(s)ds+\int_0^t\int_{\mathcal{Z}}G(s,h(s),z)\widetilde{\eta}(dz,ds)+\int_0^t\Psi(s,h(s))d W(s).
$$
For any fixed $\omega\in\Omega$, consider the deterministic
PDE:
\begin{eqnarray}\label{eq 07}
\left\{
\begin{split}
 & d M^h(t,\omega)+\mathcal{A}M^h(t,\omega)dt\\
 &\ \ \ \ \ \ \ \ \ \ \ \ \ +B\big(y_0(t,\omega),Z^h(t,\omega)M^h(t,\omega)\big)\phi_m(|y_0(t,\omega)|)g_\delta(|y_0(\omega)|_{\xi_t})dt
 = f(t)dt,\\
& M^h(0)=u_0(\omega).
\end{split}
\right.
\end{eqnarray}
According to \cite{Temam_2001}, there exists a unique $M^h(\omega)\in C([0,T],H)\cap L^2([0,T],V)$ satisfying (\ref{eq 07}), i.e.,
for every $t\in [0,T]$, as an element of $V'$,
\begin{eqnarray*}
&&\hspace{-1cm}M^h(t,\omega)\\
&=&
u_0(\omega)-\int_0^t\mathcal{A}M^h(s,\omega)ds\\
 &&-\int_0^tB\big(y_0(s,\omega),Z^h(s,\omega)+M^h(s,\omega)\big)\phi_m(|y_0(s,\omega)|)g_\delta(|y_0(\omega)|_{\xi_s})ds
 +\int_0^tf(s)ds.
\end{eqnarray*}
We see that
$\Big\{\Phi^{h}(t,\omega):=Z^h(t,\omega)+M^h(t,\omega),\ t\in[0,T],\ \omega\in\Omega\Big\}$ satisfies that
\begin{itemize}
  \item[(P1)] $\Phi^{h}\in \Upsilon_T,\ \mathbb{P}$-a.s.,
  \item[(P2)] for every $t\in [0,T]$, as an element of $V'$, $\mP\text{-}$a.s.
\begin{eqnarray}\label{eq app 1}
 \Phi^{h}(t)&=&u_0-\int_0^t\mathcal{A}\Phi^{h}(s)ds-\int_0^tB(y_0(s), \Phi^{h}(s))\phi_m(|y_0(s)|)g_\delta(|y_0|_{\xi_s})ds\nonumber
\\ && + \int_0^tf(s)ds + \int_0^{t}\int_{\mathcal{Z}}G(s,h(s),z)\widetilde{\eta}(dz,ds)+\int_0^t\Psi(s,h(s))d W(s).\nonumber
\end{eqnarray}
\end{itemize}

Applying $\rm It\hat{o}'s$ formula to $|\Phi^{h}(t)|^2$ and using (\ref{eq PB}),
one obtains
\begin{eqnarray*}
 &&|\Phi^{h}(t)|^2+2\int_0^t\|\Phi^h(s)\|^2ds\\
 &=&
 |u_0|^2+ 2\int_0^t\langle f(s), \Phi^h(s)\rangle ds
 +
 2\int_0^{t}\int_{\mathcal{Z}}\langle G(s,h(s),z), \Phi^h(s)\rangle \widetilde{\eta}(dz,ds)
 \\&&+
 2\int_0^t\langle \Psi(s,h(s)), \Phi^h(s)\rangle d W(s)
 \\&&+
 \int_0^t\int_{\mathcal{Z}}|G(s,h(s),z)|^2\eta(dz,ds)+\int_0^t\|\Psi(s,h(s))\|^2_{\mathcal{L}_2}ds.
\end{eqnarray*}
We observe that
$$
2\int_0^t\langle f(s), \Phi^h(s)\rangle ds
    \leq
        \int_0^t\|f(s)\|^2_{V'}ds
               +
        \int_0^t\|\Phi^h(s)\|^2ds,
$$
in addition,
both $\int_0^{\cdot}\int_{\mathcal{Z}}\langle G(s,h(s),z), \Phi^h(s)\rangle \widetilde{\eta}(dz,ds)$ and $\int_0^{\cdot}\langle \Psi(s,h(s)), \Phi^h(s)\rangle d W(s)$ are local martingales.
Therefore, a suitable stopping time technique (see e.g. \cite{BHZ, liuwei} ) and (\ref{eq Z01}) assure that
\begin{itemize}
  \item[(P3)]
\begin{eqnarray*}
 &&\sup_{t\in[0,T]}\mathbb{E}[|\Phi^{h}(t)|^2]
 +\mathbb{E}[\int_0^T\|\Phi^h(s)\|^2ds]\\
 &\leq&
 \mathbb{E}[|u_0|^2]+ \int_0^T\|f(s)\|^2_{V'}ds+L_3T
 +L_4T\sup_{s\in[0,T]}\mathbb{E}[|h(s)|^2]+L_5\mathbb{E}[\int_0^T\|h(s)\|^2ds]<\infty.
\end{eqnarray*}
\end{itemize}
Combining (P1)--(P3), the proof of (Claim 1) is complete.
\vskip 0.2cm

We now turn to prove (Claim 2).

Let $h_0(t)=e^{-\mathcal{A}t}u_0$, then $h_0\in\Lambda_T$.
(Claim 1) implies that we can define $h_{n+1}=\Phi^{h_n}\in\Lambda_T,\ n\geq 0$ recursively, that is, for every $t\in [0,T]$, as an element of $V'$, $\mP\text{-}$a.s.
\begin{align}\label{eq h n+1}
\begin{split}
 h_{n+1}(t)\!\!&=\!\!u_0-\int_0^t\mathcal{A}h_{n+1}(s)ds-\int_0^tB(y_0(s), h_{n+1}(s))\phi_m(|y_0(s)|)g_\delta(|y_0|_{\xi_s})ds
\\ &~~ + \int_0^tf(s)ds + \int_0^{t}\int_{\mathcal{Z}}G(s,h_n(s-),z)\widetilde{\eta}(dz,ds)+\int_0^t\Psi(s,h_n(s))d W(s).
\end{split}
\end{align}

Next we will estimate $h_{n+1}(t)-h_n(t)$.

By (\ref{eq h n+1}), (\ref{eq PB}), and $\rm It\hat{o}'s$ formula, we get
\begin{eqnarray*}
 &&|h_{n+1}(t)-h_n(t)|^2+2\int_0^t\|h_{n+1}(s)-h_n(s)\|^2ds\\
 &=&
 2\int_0^{t}\int_{\mathcal{Z}}\langle G(s,h_n(s-),z)-G(s,h_{n-1}(s-),z), h_{n+1}(s-)-h_n(s-)\rangle \widetilde{\eta}(dz,ds)
 \\&&+
 2\int_0^t\langle \Psi(s,h_n(s))-\Psi(s,h_{n-1}(s)), h_{n+1}(s)-h_n(s))\rangle d W(s)
 \\&&+
 \int_0^t\int_{\mathcal{Z}}|G(s,h_n(s-),z)-G(s,h_{n-1}(s-),z)|^2\eta(dz,ds)
 \\&& +\int_0^t\|\Psi(s,h_n(s))-\Psi(s,h_{n-1}(s))\|^2_{\mathcal{L}_2}ds\\
 &=:&\sum_{i=1}^4I_i(t).
\end{eqnarray*}
 Moreover, $I_1$ and $I_2$ are local martingales.
 This and
 using a suitable stopping time technique hence
 demonstrates that
\begin{eqnarray*}
 &&\sup_{t\in[0,T]}\mathbb{E}[|h_{n+1}(t)-h_n(t)|^2]+2\mathbb{E}[\int_0^T\|h_{n+1}(s)-h_n(s)\|^2ds]\\
 &\leq&
 \mathbb{E}[\int_0^T\int_{\mathcal{Z}}|G(s,h_n(s),z)-G(s,h_{n-1}(s),z)|^2\nu(dz)ds]\\
 &&+\mathbb{E}[\int_0^T\|\Psi(s,h_n(s))-\Psi(s,h_{n-1}(s))\|^2_{\mathcal{L}_2}ds]\\
 &\leq&
 L_1\mathbb{E}[\int_0^T|h_n(s)-h_{n-1}(s)|^2ds]
 +
 L_2\mathbb{E}[\int_0^T\|h_n(s)-h_{n-1}(s)\|^2ds]\\
 &\leq&
 L_1T\sup_{s\in[0,T]}\mathbb{E}[|h_n(s)-h_{n-1}(s)|^2]
 +
 L_2\mathbb{E}[\int_0^T\|h_n(s)-h_{n-1}(s)\|^2ds].
\end{eqnarray*}
Here Condition \ref{con G} has been used to get the second inequality. Multiplying $\frac{1}{2}$ at both sides of the above inequality, we obtain
\begin{eqnarray*}
 &&\sup_{t\in[0,T]}\mathbb{E}[\frac{1}{2}|h_{n+1}(t)-h_n(t)|^2]+\mathbb{E}[\int_0^T\|h_{n+1}(s)-h_n(s)\|^2ds]\\
 &\leq&
 (L_1 T\vee \frac{L_2}{2})
 \Big(
 \sup_{s\in[0,T]}\mathbb{E}[\frac{1}{2}|h_n(s)-h_{n-1}(s)|^2]
 +
\mathbb{E}[\int_0^T\|h_n(s)-h_{n-1}(s)\|^2ds]
\Big).
\end{eqnarray*}

Choosing $T_0>0$ such that $L_1T_0<1$ and noticing that $L_2<2$ (see Condition \ref{con G}).
Thus,
there exists {an} $H$-valued process $\Theta$ such that 
{$\Theta$} has a $\mathbb{F}$-progressively measurable version, denoted by $\widetilde{\Theta}$,
\begin{eqnarray}\label{eq BJ 01}
\lim_{n\rightarrow\infty}\sup_{t\in[0,{T_0}]}\mathbb{E}[|h_n(t)-\Theta(t)|^2]=0,
\end{eqnarray}
and
\begin{eqnarray}\label{eq BJ 02}
\lim_{n\rightarrow\infty}\mathbb{E}[\int_0^{T_0}\|h_n(t)-\Theta(t)\|^2dt]=0.
\end{eqnarray}
Note that (\ref{eq BJ 01}), (\ref{eq BJ 02}) and Condition \ref{con G}
assure that
\begin{itemize}
  \item[(Q1)]
    \begin{eqnarray*}
        &&\mathbb{E}[\int_0^{T_0}\int_{\mathcal{Z}}|G(s,h_n(s),z)-G(s,\Theta(s),z)|^2\nu(dz)ds]\\
       && +
        \mathbb{E}[\int_0^{T_0}\|\Psi(s,h_n(s))-\Psi(s,\Theta(s))\|^2_{\mathcal{L}_2}ds]\\
     &\leq&
        L_1\mathbb{E}[\int_0^{T_0}|h_n(s)-\Theta(s)|^2ds]
            +
        L_2\mathbb{E}[\int_0^{T_0}\|h_n(s)-\Theta(s)\|^2ds]\\
     &\leq&
      L_1{T_0}\sup_{s\in[0,T_0]}\mathbb{E}[|h_n(s)-\Theta(s)|^2]
            +
        L_2\mathbb{E}[\int_0^{T_0}\|h_n(s)-\Theta(s)\|^2ds]\\
     &&\rightarrow 0,\ \ \text{as }n\rightarrow\infty.
    \end{eqnarray*}
  \item[(Q2)] $\mathbb{E}[\int_0^{T_0}\|\mathcal{A}h_n(s)-\mathcal{A}\Theta(s)\|^2_{V'}ds]
               =
               \mathbb{E}[\int_0^{T_0}\|h_n(s)-\Theta(s)\|^2ds]\rightarrow 0,\ \ \text{as }n\rightarrow\infty.
               $
  \item[(Q3)] For any $t\in[0,T_0]$ and $e\in V$, by Condition \ref{p-2},
  \begin{eqnarray*}
     &&\mathbb{E}[|\int_0^t\langle B(y_0(s), h_n(s))\phi_m(|y_0(s)|)g_\delta(|y_0|_{\xi_s})
                                \\
&&\ \ \ \  -
                                B(y_0(s), \Theta(s))\phi_m(|y_0(s)|)g_\delta(|y_0|_{\xi_s}),
                                e
                        \rangle ds|]\\
     &\leq&
     \mathbb{E}[\int_0^{T_0}|y_0(s)|^{\frac{1}{2}}\|y_0(s)\|^{\frac{1}{2}}|e|^{\frac{1}{2}}\|e\|^{\frac{1}{2}}\phi_m(|y_0(s)|)g_\delta(|y_0|_{\xi_s}) \|h_n(s)-\Theta(s)\|ds]\\
     &\leq&
     (m+1)^{\frac{1}{2}}|e|^{\frac{1}{2}}\|e\|^{\frac{1}{2}}T_0^{\frac{1}{4}}
     \Big(\mathbb{E}[\int_0^{T_0} \|h_n(s)-\Theta(s)\|^2ds]\Big)^{\frac{1}{2}}\\
     &&\ \  \cdot
     \Big(\mathbb{E}[\int_0^{T_0} \|y_0(s)\|^2g_\delta(|y_0|_{\xi_s})ds]\Big)^{\frac{1}{4}}\\
     &\leq&
     \Big(2\delta(m+1)|e|\|e\|\Big)^{\frac{1}{2}}T_0^{\frac{1}{4}}
     \Big(\mathbb{E}[\int_0^{T_0} \|h_n(s)-\Theta(s)\|^2ds]\Big)^{\frac{1}{2}}\rightarrow 0,\ \ \text{as }n\rightarrow\infty.
  \end{eqnarray*}

  \item[(Q4)] For any $t\in[0,T_0]$, there exists a subsequence $n_k\uparrow\infty$ such that
  $$
  \lim_{k\rightarrow\infty}|h_{n_k}(t)-\Theta(t)|=0,\ \ \mathbb{P}\text{-}a.s..
  $$

\end{itemize}

By the definition of $h_{n+1}:=\Phi^{h_n}$(see (P2)),  for every $t\in [0,T_0]$, $\mP\text{-}$a.s. for any $e\in V$,
\begin{eqnarray*}
  && \langle h_{n+1}(t),e\rangle\\
&=&
   \langle u_0,e\rangle
      -
   \int_0^t\langle\mathcal{A}h_{n+1}(s),e\rangle ds
      +
   \int_0^t\langle B(y_0(s), h_{n+1}(s))\phi_m(|y_0(s)|)g_\delta(|y_0|_{\xi_s}),e\rangle ds\nonumber
\\ \!\!&& +\!\!
   \int_0^t\langle f(s),e\rangle ds
     +
   \int_0^{t}\int_{\mathcal{Z}}\langle G(s,h_n(s-),z),e\rangle \widetilde{\eta}(dz,ds)
     +
   \int_0^t\langle \Psi(s,h_n(s)),e\rangle d W(s).\nonumber
\end{eqnarray*}
Applying (Q1)-(Q4) and taking limits in the above equation (choosing a subsequence if necessary), we obtain
\begin{eqnarray*}
 &&  \langle \Theta(t),e\rangle\\
&=&
   \langle u_0,e\rangle
      -
   \int_0^t\langle\mathcal{A}\Theta(s),e\rangle ds
      +
   \int_0^t\langle B(y_0(s), \Theta(s))\phi_m(|y_0(s)|)g_\delta(|y_0|_{\xi_s}),e\rangle ds\nonumber
\\ && +
   \int_0^t\langle f(s),e\rangle ds
     +
   \int_0^{t}\int_{\mathcal{Z}}\langle G(s,\widetilde{\Theta}(s),z),e\rangle \widetilde{\eta}(dz,ds)
     +
   \int_0^t\langle \Psi(s,\widetilde{\Theta}(s)),e\rangle d W(s).\nonumber
\end{eqnarray*}

Applying $\rm It\hat{o}'s$ formula to $|\Theta(t)|^2$, it holds that
\begin{eqnarray*}
  && |\Theta(t)|^2+2\int_0^t\|\Theta(s)\|^2ds\\
&=&
   |u_0|^2
 +
   2\int_0^t\langle f(s),\Theta(s)\rangle ds
     +
   2\int_0^{t}\int_{\mathcal{Z}}\langle G(s,\widetilde{\Theta}(s),z),\widetilde{\Theta}(s)\rangle \widetilde{\eta}(dz,ds)
   \\
 &&  +
   2\int_0^t\langle \Psi(s,\widetilde{\Theta}(s)),\widetilde{\Theta}(s)\rangle d W(s)+
    \int_0^{t}\int_{\mathcal{Z}}|G(s,\widetilde{\Theta}(s),z)|^2\eta(dz,ds)
        +
    \int_0^{t}\|\Psi(s,\widetilde{\Theta}(s))\|^2_{\mathcal{L}_2}ds.
\end{eqnarray*}
Since $\int_0^{\cdot}\int_{\mathcal{Z}}\langle G(s,\widetilde{\Theta}(s),z),\widetilde{\Theta}(s)\rangle \widetilde{\eta}(dz,ds)$ and $2\int_0^{\cdot}\langle \Psi(s,\widetilde{\Theta}(s)),\widetilde{\Theta}(s)\rangle d W(s)$ are local martingales,  using a suitable stopping time technique again, the above inequality implies that
\begin{eqnarray*}
  && \mathbb{E}[|\Theta(t)|^2]+2\mathbb{E}[\int_0^t\|\Theta(s)\|^2ds]\\
&\leq&
   \mathbb{E}[|u_0|^2]
 +
   2\mathbb{E}[\int_0^t\langle f(s),\Theta(s)\rangle ds]
 +
    \mathbb{E}[\int_0^{t}\int_{\mathcal{Z}}|G(s,\widetilde{\Theta}(s),z)|^2\nu(dz)ds]
       \\
       && +
    \mathbb{E}[\int_0^{t}\|\Psi(s,\widetilde{\Theta}(s))\|^2_{\mathcal{L}_2}ds]\\
 &\leq&
    \mathbb{E}[|u_0|^2]
 +
   \epsilon \mathbb{E}[\int_0^t\|\Theta(s)\|^2ds]
   +
   \epsilon^{-1}\int_0^T\|f(s)\|^2_{V'}ds\\
    &&+
    L_3t
    +
    L_4\mathbb{E}[\int_0^t|\Theta(s)|^2ds]
        +
    L_5\mathbb{E}[\int_0^t\|\Theta(s)\|^2ds].
\end{eqnarray*}
Choosing $\epsilon=\frac{2-L_5}{2}$, the above inequality shows that
\begin{eqnarray*}
  && \mathbb{E}[|\Theta(t)|^2]+\frac{2-L_5}{2}\mathbb{E}[\int_0^t\|\Theta(s)\|^2ds]\\
 &\leq&
    \mathbb{E}[|u_0|^2]
   +
   \frac{2}{2-L_5}\int_0^T\|f(s)\|^2_{V'}ds
    +
    L_3t
    +
    L_4\int_0^t\mathbb{E}[|\Theta(s)|^2]ds.
\end{eqnarray*}
Gronwall's lemma ensures that
\begin{eqnarray}\label{eq BJ 04}
  && \sup_{t\in[0,T_0]}\mathbb{E}[|\Theta(t)|^2]
  +\frac{2-L_2}{2}\mathbb{E}[\int_0^{T_0}\|\Theta(s)\|^2ds]\nonumber\\
 &\leq&
    \Big(\mathbb{E}[|u_0|^2]
   +
   \frac{2}{2-L_5}\int_0^{T_0}\|f(s)\|^2_{V'}ds
   +L_3T_0\Big)\exp(L_1T_0).
\end{eqnarray}
(Claim 1) demonstrates that $\Theta\in \Lambda_{T_0}$ and $\Theta$ is a solution of (\ref{eq BJ00}) on $[0,T_0]$. By a standard argument, for any $T>0$, there exists a solution of (\ref{eq BJ00}) on $[0,T]$. The uniqueness is standard, and thus we omit it here. (\ref{eq BJ 04}) implies that
(\ref{eq BJ 03}) holds.
Hence, the statements in
(Claim 2) are proved, and
the proof of Lemma \ref{Lemma 1} is thus complete.
\hfill$\Box$\\

Set 
{ $y_0(t):=0$}.  By Lemma \ref{Lemma 1},
{for any $m\in \mN$ and $\delta>0,$  we can define the sequence $\{y_n\}_{n=1}^\infty$ recursively} by $y_{n+1}:=\Theta^{y_n}$, which satisfies the following equation
\begin{eqnarray}\label{XZ-1}
\left\{
\begin{split}
& d y_{n+1}(t) +\mathcal{A}y_{n+1}(t)dt
+B(y_{n}(t),y_{n+1}(t))\phi_m(|y_n(t)|)
g_\delta(|y_n|_{\xi_t})dt
\\ & = f(t)dt + \int_{\mathcal{Z}}G(t,y_{n+1}(t-),z)\widetilde{\eta}(dz,dt)+\Psi(t,y_{n+1}(t))d W(t),
\\
& y_{n+1}(0)=u_0.
\end{split}
\right.
\end{eqnarray}
Moreover,
\begin{eqnarray}\label{eq BJ 05}
 \sup_{t\in[0,T]}\mathbb{E}[|y_n(t)|^2]
 +\mathbb{E}[\int_0^{T}\|y_n(s)\|^2ds]
 \leq
 C\Big(\mathbb{E}[|u_0|^2],
   \int_0^T\|f(s)\|^2_{V'}ds, T\Big).
\end{eqnarray}
Here $C\Big(\mathbb{E}[|u_0|^2],
   \int_0^T\|f(s)\|^2_{V'}ds, T\Big)$ is independent of 
   {$m,\delta,n$}.

Note that
\begin{eqnarray}\label{eq y1}
\left\{
\begin{split}
& d y_{1}(t) +\mathcal{A}y_{1}(t)dt
= f(t)dt + \int_{\mathcal{Z}}G(t,y_{1}(t-),z)\widetilde{\eta}(dz,dt)+\Psi(t,y_{1}(t))d W(t),
\\
& y_{1}(0)=u_0.
\end{split}
\right.
\end{eqnarray}

\vskip 0.2cm

We will prove that there exist $\delta_0>0$ and $T_0>0$ such that
\begin{eqnarray}\label{eq LIp}
\sum_{n=2}^\infty\Big(\mathbb{E}[\int_0^{T_0}\|y_{n+1}(t)-y_n(t)\|^2ds]\Big)^{1/2}
+
\sum_{n=2}^\infty\mathbb{E}[\sup_{t\in[0,T_0]}|y_{n+1}(t)-y_n(t)|]<\infty.
\end{eqnarray}
Here $y_{n+1}$ is the solution of (\ref{XZ-1}) with $\delta$ replaced by $\delta_0$.
The proof of this claim needs Lemmas~\ref{lemma 2} and \ref{lemma 3} below.
To reduce the proof, we need some notations.
\vskip 0.2cm
Let
\begin{eqnarray}
\label{1-1}
\begin{split}
 I_n(t)&=
 \Big\langle B(y_{n}(t),y_{n+1}(t))\phi_m(|y_n(t)|)g_\delta(|y_n|_{\xi_t}) \\&\quad\quad  -B(y_{n-1}(t),y_n(t))\phi_m(|y_{n-1}(t)|)g_\delta(|y_{n-1}|_{\xi_t}),
  y_{n+1}(t)-y_n(t)
\Big\rangle.
\end{split}
\end{eqnarray}
Let  $I_{[0,l]}:(-\infty, \infty)\rightarrow \{0,1\}$ be an indicator function defined by
\begin{eqnarray*}
I_{[0,l]}(x)=\left\{\begin{split}
 & 1, \ \ \text{ if } x\in [0,l],
 \\
 & 0, \ \ \text{ else} .
\end{split}
\right.
\end{eqnarray*}
For any $n\in \mN,\eps,p,t>0,$ we set
\begin{eqnarray}\label{eq Xi}
\Xi_n(t):=\|y_{n-1}(t)\|^{2}I_{[0,3\delta]}(|y_{n-1}|_{\xi_t})+\|y_n(t)\|^2 I_{[0,3\delta]}(|y_n|_{\xi_t}),
\end{eqnarray}
and
{\begin{eqnarray}\label{eq Sn}
\$_n(t):=\eps^{-3}\Big(1+(m+2)^2+\delta^{-1}(m+2)^2+(m+1)^2\delta^{-4p}+(m+1)^2 \eps^3 \delta^{-4p}\Big)\Xi_n(t).
\end{eqnarray}
}
\begin{lemma}\label{lemma 2}
{For any $\eps, p,t>0$ and $ n\geq 1$},  the following inequality holds:
\begin{eqnarray}\label{Eq In}
\begin{split}
  I_n(t) & \leq    7\eps   \|y_{n+1}(t)-y_n(t)\|^{2}+\Big(2\eps+\eps^{-1/2}\delta^{2p}\Big)   \|y_n(t)-y_{n-1}(t)\|^{2}
\\  & +C  \Big(\frac{\eps  }{\delta^{3/2}}+\eps^{-1/2}\delta^{2p-2}+\eps \delta^{2(p-1)}\Big)
|y_n-y_{n-1}|_{\xi_t}^{2} \Xi_n(t)
\\  &+3\eps |y_n(t)-y_{n-1}(t)|^{2}
      \Xi_n(t)
+{ C \$_n (t)}|y_{n+1}(t)-y_n(t)|^{2},
\end{split}
\end{eqnarray}
{where   $C$  is a constant independent of $\eps,\delta, n,p,t$.}
\end{lemma}

\noindent\emph{Proof of Lemma \ref{lemma 2}.}\quad
We will prove this lemma in the  following   four different cases.



\begin{enumerate}
  \item [(\uppercase\expandafter{\romannumeral1})]:  $|y_n|_{\xi_t}\leq 3\delta$ and  $|y_{n-1}|_{\xi_t}\leq 3\delta$,
  \item [(\uppercase\expandafter{\romannumeral2})]:  $|y_n|_{\xi_t}\leq 3\delta$ and  $|y_{n-1}|_{\xi_t}>3\delta$,
  \item [(\uppercase\expandafter{\romannumeral3})]: $|y_n|_{\xi_t}> 3\delta$ and  $|y_{n-1}|_{\xi_t}\leq 3\delta$,
  \item [(\uppercase\expandafter{\romannumeral4})]: $|y_n|_{\xi_t}> 3\delta$ and  $|y_{n-1}|_{\xi_t}> 3\delta$.
\end{enumerate}

\textbf{(\uppercase\expandafter{\romannumeral1}):}  \textbf{$|y_n|_{\xi_t}\leq 3\delta$ and  $|y_{n-1}|_{\xi_t}\leq 3\delta$.}

We need to divide this case \textbf{(\uppercase\expandafter{\romannumeral1})}   into the following four  subcases.
\begin{enumerate}
  \item [(\uppercase\expandafter{\romannumeral1}-1)]:  $|y_n(t)|\leq m+2$ and  $|y_{n-1}(t)|\leq m+2$,
  \item [(\uppercase\expandafter{\romannumeral1-2})]:  $|y_n(t)|\leq  m+2$ and  $|y_{n-1}(t)|> m+2$,
  \item [(\uppercase\expandafter{\romannumeral1-3})]:  $|y_n(t)|> m+2$ and  $|y_{n-1}(t)|\leq m+2$,
  \item [(\uppercase\expandafter{\romannumeral1-4})]:  $|y_n(t)|> m+2$ and  $|y_{n-1}(t)|> m+2$.
\end{enumerate}

\emph{(\uppercase\expandafter{\romannumeral1}-1):}
  $|y_n(t)|\leq m+2$ and  $|y_{n-1}(t)|\leq m+2$.

Observing  (\ref{p-3}) and (\ref{eq PB}),  we have
\begin{eqnarray}
\nonumber
&& I_n(t)
\\ \nonumber  && = \langle B(y_{n}(t),y_{n+1}(t))\phi_m(|y_n(t)|)g_\delta(|y_n|_{\xi_t}) -B(y_{n}(t),y_n(t))\phi_m(|y_{n}(t)|)g_\delta(|y_{n} |_{\xi_t}),
\\  \nonumber  &&\quad \quad  y_{n+1}(t)-y_n(t) \rangle
\\  \nonumber  && +\langle B(y_{n}(t),y_{n}(t))\phi_m(|y_n(t)|)g_\delta(|y_n|_{\xi_t}) -B(y_{n-1}(t),y_n(t))\phi_m(|y_{n}(t)|)g_\delta(|y_{n} |_{\xi_t}),
\\ \nonumber  &&\quad \quad  y_{n+1}(t)-y_n(t) \rangle
\\ \nonumber  &&+ \langle B(y_{n-1}(t),y_n(t))\phi_m(|y_{n}(t)|)g_\delta(|y_{n} |_{\xi_t})-B(y_{n-1}(t),y_n(t))\phi_m(|y_{n-1}(t)|)g_\delta(|y_{n-1} |_{\xi_t}),
\\ \nonumber  &&\quad \quad  y_{n+1}(t)-y_n(t) \rangle
\\ \label{p-6}&&:=0+J_1(t)+J_2(t).
\end{eqnarray}
We will give an estimate of $J_1(t)$ and  $J_2(t)$ respectively.
By (\ref{p-3}) (\ref{p-4})(\ref{p-5}) and the Young inequality,  for any $\eps>0$, we have
\begin{eqnarray}
\nonumber   && J_1(t)\\
&\leq &\Big| \langle B(y_{n}(t)-y_{n-1}(t),y_{n}(t)) ,  y_{n+1}(t)-y_n(t) \rangle\Big|
  \\  \nonumber  &\leq &  C |y_n(t)-y_{n-1}(t)|_Q \cdot \|y_n(t)\|\cdot  |y_{n+1}(t)-y_n(t)|_Q
  \\ \nonumber  &\leq & C|y_n(t)-y_{n-1}(t)|^{1/2}\|y_n(t)-y_{n-1}(t)\|^{1/2} \cdot \|y_{n+1}(t)-y_{n}\|^{1/2}\cdot |y_{n+1}(t)-y_{n}|^{1/2}\cdot \|y_n(t)\|
  \\ \nonumber  &\leq & \eps \|y_{n+1}(t)-y_n(t)\| \cdot \|y_n(t)-y_{n-1}(t)\|+C\eps^{-1} \|y_n(t)\|^2  |y_{n+1}(t)-y_n(t)|\cdot |y_n(t)-y_{n-1}(t)|
  \\  \nonumber  &\leq & \eps \|y_{n+1}(t)-y_n(t)\|^2+\eps  \|y_n(t)-y_{n-1}(t)\|^2+C\eps^{-3}\|y_n(t)\|^2  |y_{n+1}(t)-y_n(t)|^2
  \\
   &&+\eps  \|y_n(t)\|^2 |y_n(t)-y_{n-1}(t)|^2,
    \label{p-7}
\end{eqnarray}
\begin{eqnarray}
 \nonumber   J_2(t)&=&\Big\langle B(y_{n-1}(t),y_{n}(t))\big[\phi_m(|y_{n}(t)|)g_\delta(|y_{n} |_{\xi_t})-\phi_m(|y_{n-1}(t)|)g_\delta(|y_{n-1} |_{\xi_t})\big],
  \\ \nonumber   && \quad  y_{n+1}(t)-y_n(t) \Big\rangle
 \\  \nonumber  &=& \Big\langle B(y_{n-1}(t),y_{n}(t))\big[\phi_m(|y_{n}(t)|)g_\delta(|y_{n} |_{\xi_t})-\phi_m(|y_{n-1}(t)|)g_\delta(|y_{n} |_{\xi_t})\big],
  \\  \nonumber  && \quad  y_{n+1}(t)-y_n(t) \Big\rangle
  \\  \nonumber  &&+
  \Big\langle B(y_{n-1}(t),y_{n}(t))\big[\phi_m(|y_{n-1}(t)|)g_\delta(|y_{n} |_{\xi_t})-\phi_m(|y_{n-1}(t)|)g_\delta(|y_{n-1} |_{\xi_t})\big],
  \\ \nonumber   && \quad  y_{n+1}(t)-y_n(t) \Big\rangle
 \\ \label{p-8} :&=&J_{2,1}(t)+J_{2,2}(t).
 \end{eqnarray}
 (\ref{p-4}), (\ref{p-5}), the Lipchiz property of $\phi_m$ and $g_\delta$, and the Young inequality prove that for any $\eps>0,$
 \begin{eqnarray}\label{p-9}
    J_{2,1}(t)
 &\leq & C|y_n(t)-y_{n-1}(t)| |y_{n-1}(t)|_Q \cdot \|y_n(t)\| \cdot
    |y_{n+1}(t)-y_n(t)|_Q  \nonumber
   \\  \nonumber
 &\leq & C|y_n(t)-y_{n-1}(t)|\|y_{n+1}(t)-y_n(t)\|^{1/2}\cdot |y_{n+1}(t)-y_n(t)|^{1/2} \\  \nonumber
 &&\cdot \|y_{n-1}(t)\|^{1/2}|y_{n-1}(t)|^{1/2} \|y_n(t)\|
 \\  \nonumber  &\leq & C\eps^{-3}|y_{n+1}(t)-y_n(t)|^{2}\|y_{n-1}(t)\|^{2}|y_{n-1}(t)|^{2}
  \\  \nonumber
  &&+\eps|y_n(t)-y_{n-1}(t)|^{4/3}\|y_n(t)\|^{4/3}\|y_{n+1}(t)-y_n(t)\|^{2/3}
 \\  \nonumber  &\leq & C\eps^{-3}(m+2)^2|y_{n+1}(t)-y_n(t)|^{2}\|y_{n-1}(t)\|^{2}
 +\eps|y_n(t)-y_{n-1}(t)|^{2}\|y_n(t)\|^{2} \\
 &&+\eps \|y_{n+1}(t)-y_n(t)\|^{2}.
\end{eqnarray}
and
\begin{eqnarray}
 J_{2,2}(t)
 &\leq &  C\frac{1}{\delta} |y_n-y_{n-1}|_{\xi_t} |y_{n-1}(t)|_Q \cdot \|y_n(t)\| \cdot
    |y_{n+1}(t)-y_n(t)|_Q \nonumber
 \\ \nonumber   &\leq & C\frac{1}{\delta} |y_n-y_{n-1}|_{\xi_t}\|y_{n+1}(t)-y_n(t)\|^{1/2}\cdot |y_{n+1}(t)-y_n(t)|^{1/2}
 \\ &&\cdot \|y_{n-1}(t)\|^{1/2}|y_{n-1}(t)|^{1/2} \|y_n(t)\|
 \\  \nonumber  &\leq & C\eps^{-3}\Big(\frac{1}{\delta}\Big)^{}|y_{n+1}(t)-y_n(t)|^{2}\|y_{n-1}(t)\|^{2}|y_{n-1}(t)|^{2}
  \\  \nonumber
  &&
 +\eps \Big(\frac{1}{\delta}\Big)^{} |y_n-y_{n-1}|_{\xi_t}^{4/3}\|y_n(t)\|^{4/3}\|y_{n+1}(t)-y_n(t)\|^{2/3}
 \\ \nonumber  &\leq &  C\eps^{-3}(m+2)^2 \Big(\frac{1}{\delta}\Big) |y_{n+1}(t)-y_n(t)|^{2}\|y_{n-1}(t)\|^{2}
 +\eps \Big(\frac{1}{\delta}\Big)^{3/2}|y_n-y_{n-1}|_{\xi_t}^{2}\|y_n(t)\|^{2}
 \\ \label{p-10}&&\quad +\eps \|y_{n+1}(t)-y_n(t)\|^{2}.
\end{eqnarray}
Combining 
(\ref{p-7})--\eqref{p-10} with (\ref{p-6}), and $|y_n|_{\xi_t}\leq 3\delta$ and  $|y_{n-1}|_{\xi_t}\leq 3\delta$, one yields that for any $\eps>0$, the following inequality holds for this subcase,
\begin{eqnarray}\label{p-15}
 I_n(t)
 \nonumber  & \leq &   3\eps   \|y_{n+1}(t)-y_n(t)\|^{2}+\eps   \|y_n(t)-y_{n-1}(t)\|^{2}
\\ \nonumber && + \frac{\eps  }{\delta^{3/2}}\|y_n(t)\|^2  |y_n-y_{n-1}|_{\xi_t}^{2} I_{[0,3\delta]}(|y_n|_{\xi_t})
\\ \nonumber&&+2\eps  \|y_n(t)\|^2 |y_n(t)-y_{n-1}(t)|^2 I_{[0,3\delta]}(|y_n|_{\xi_t})
\\  \nonumber &&+C\eps^{-3}\big(1+(m+2)^2+\delta^{-1}(m+2)^2\big)|y_{n+1}(t)-y_n(t)|^{2}\Xi_n(t).
\end{eqnarray}

\emph{(\uppercase\expandafter{\romannumeral1}-2): $|y_n(t)|\leq  m+2$ and  $|y_{n-1}(t)|> m+2$.}

For this subcase, according to the defintion of $\phi_m$, it holds that
\begin{eqnarray*}
&& I_n(t)\cdot I_{\{|y_n|_{\xi_t}\leq 3\delta\}}\cdot I_{\{|y_{n-1}|_{\xi_t}\leq 3\delta\}}\cdot I_{\{|y_n(t)|\leq  m+2\}}\cdot I_{\{|y_{n-1}(t)|> m+2\}}\\
&=&
  I_n(t)\cdot I_{\{|y_n|_{\xi_t}\leq 3\delta\}}\cdot I_{\{|y_{n-1}|_{\xi_t}\leq 3\delta\}}\cdot I_{\{|y_n(t)|\leq  m+1\}}\cdot I_{\{|y_{n-1}(t)|> m+2\}}.
\end{eqnarray*}
For any $t$ such that $I_{\{|y_n|_{\xi_t}\leq 3\delta\}}\cdot I_{\{|y_{n-1}|_{\xi_t}\leq 3\delta\}}\cdot I_{\{|y_n(t)|\leq  m+1\}}\cdot I_{\{|y_{n-1}(t)|> m+2\}}=1$, we have
\begin{eqnarray}\label{p-12}
|y_n(t)-y_{n-1}(t)|\geq  1.
\end{eqnarray}
This and (\ref{p-3}), (\ref{p-4})-(\ref{eq PB}) show that 
\begin{eqnarray}
I_n(t)
 \nonumber &=&\langle B(y_{n}(t),y_{n+1}(t))\phi_m(|y_n(t)|)g_\delta(|y_n|_{\xi_t}) ,  y_{n+1}(t)-y_n(t) \rangle
\\  \nonumber &\leq& \big| \langle B(y_{n}(t),y_{n+1}(t)) ,  -y_n(t) \rangle\big|
\\  \nonumber &=&\big| \langle B(y_{n}(t),y_{n+1}(t)-y_n(t)),  -y_n(t) \rangle\big|
\\ \nonumber  &=& \big| \langle B(y_{n}(t),y_n(t) ), y_{n+1}(t)-y_n(t)\rangle\big|
\\  \nonumber &\leq&  C |y_n(t)|_Q \cdot \|y_n(t)\| \cdot |y_{n+1}(t)-y_n(t)|_Q
\\ \nonumber &\leq&  C\|y_{n+1}(t)-y_n(t)\|^{1/2} |y_{n+1}(t)-y_n(t)|^{1/2}\cdot \|y_n(t)\| \cdot |y_n(t)|^{1/2}
\|y_n(t)\|^{1/2}
\end{eqnarray}
(\ref{p-12}) and Young's inequality therefore assure that for any $\epsilon>0$,
\begin{eqnarray}
I_n(t)
&\leq&  C\|y_{n+1}(t)-y_n(t)\|^{1/2} |y_{n+1}(t)-y_n(t)|^{1/2}\cdot \|y_n(t)\| \nonumber\\
&&\cdot |y_n(t)|^{1/2}
\|y_n(t)\|^{1/2}\cdot  |y_n(t)-y_{n-1}(t)|\nonumber\\
 &\leq&  C(m+1)^{1/2} \|y_{n+1}(t)-y_n(t)\|^{1/2} |y_{n+1}(t)-y_n(t)|^{1/2}\cdot \|y_n(t)\| \nonumber\\  \nonumber
 &&\cdot
\|y_n(t)\|^{1/2}\cdot  |y_n(t)-y_{n-1}(t)|
\\ \nonumber
&\leq &C(m+1)^2\eps^{-3} \|y_n(t)\|^2  |y_{n+1}(t)-y_n(t)|^2
\\  \nonumber
&&+\eps  \|y_{n+1}(t)-y_n(t)\|^{2/3}\|y_n(t)\|^{4/3} |y_n(t)-y_{n-1}(t)|^{4/3}
\\  \nonumber
&\leq & C(m+1)^2\eps^{-3} \|y_n(t)\|^2  |y_{n+1}(t)-y_n(t)|^2I_{[0,3\delta]}(|y_n|_{\xi_t})
+\eps  \|y_{n+1}(t)-y_n(t)\|^{2}
\\  \label{p-16} && \quad +\eps \|y_n(t)\|^{2} |y_n(t)-y_{n-1}(t)|^{2}I_{[0,3\delta]}(|y_n|_{\xi_t}).
\end{eqnarray}
To get the last inequality, we have used the fact that for this subcase $|y_n|_{\xi_t}\leq 3\delta$ and  $|y_{n-1}|_{\xi_t}\leq 3\delta$.

\vskip 0.2cm
\emph{(\uppercase\expandafter{\romannumeral1}-3):   $|y_n(t)|> m+2$ and  $|y_{n-1}(t)|\leq m+2$.}

In  this subcase, by the defintion of $\phi_m$, we obtain
\begin{eqnarray*}
&& I_n(t)\cdot I_{\{|y_n|_{\xi_t}\leq 3\delta\}}\cdot I_{\{|y_{n-1}|_{\xi_t}\leq 3\delta\}}\cdot I_{\{|y_n(t)|>  m+2\}}\cdot I_{\{|y_{n-1}(t)|\leq m+2\}}\\
&=&
  I_n(t)\cdot I_{\{|y_n|_{\xi_t}\leq 3\delta\}}\cdot I_{\{|y_{n-1}|_{\xi_t}\leq 3\delta\}}\cdot I_{\{|y_n(t)|>  m+2\}}\cdot I_{\{|y_{n-1}(t)|\leq m+1\}}.
\end{eqnarray*}
For any $t$ such that $I_{\{|y_n|_{\xi_t}\leq 3\delta\}}\cdot I_{\{|y_{n-1}|_{\xi_t}\leq 3\delta\}}\cdot  I_{\{|y_n(t)|>  m+2\}}\cdot I_{\{|y_{n-1}(t)|\leq m+1\}}=1$, we have
\begin{eqnarray}\label{p-13}
|y_n(t)-y_{n-1}(t)|\geq  1.
\end{eqnarray}
Moreover, (\ref{p-4}), (\ref{p-5}) ensure that for any $\eps>0$,
\begin{eqnarray}
\nonumber |I_n(t)|
 \nonumber &=&|\langle B(y_{n-1}(t),y_{n}(t))\phi_m(|y_{n-1}(t)|) g_\delta(|y_{n-1}|_{\xi_t}),  y_{n+1}(t)-y_n(t) \rangle|
\\  \nonumber &\leq &C |y_{n-1}(t)|_Q \cdot  \|y_n(t)\|\cdot | y_{n+1}(t)-y_n(t)|_Q
\\ \nonumber &\leq & C \|y_{n+1}(t)-y_n(t)\|^{1/2}|y_{n+1}(t)-y_n(t)|^{1/2}\cdot\|y_n(t)\|\cdot
\|y_{n-1}(t)\|^{1/2}|y_{n-1}(t)|^{1/2}
\\ \nonumber&\leq & C(m+1)^{1/2}  \|y_{n+1}(t)-y_n(t)\|^{1/2}|y_{n+1}(t)-y_n(t)|^{1/2}\cdot\|y_n(t)\|\cdot
\|y_{n-1}(t)\|^{1/2}
\\ \nonumber &\leq & C\eps^{-3}(m+1)^2|y_{n+1}(t)-y_n(t)|^{2} \|y_{n-1}(t)\|^{2}+\eps \|y_n(t)\|^{4/3}
\|y_{n+1}(t)-y_n(t)\|^{2/3}
\\ \nonumber &\leq & C\eps^{-3}(m+1)^2|y_{n+1}(t)-y_n(t)|^{2} \|y_{n-1}(t)\|^{2}
+\eps \|y_n(t)\|^{2}
+\eps \|y_{n+1}(t)-y_n(t)\|^{2}
\\ \nonumber &\leq &   C\eps^{-3}(m+1)^2|y_{n+1}(t)-y_n(t)|^{2} \|y_{n-1}(t)\|^{2}
+\eps \|y_n(t)-y_{n-1}(t)\|^{2}+\eps \|y_{n-1}(t)\|^2
\\ \nonumber
&&+\eps \|y_{n+1}(t)-y_n(t)\|^{2}.
\end{eqnarray}
Combining this with (\ref{p-13}) imply that
\begin{eqnarray}
\nonumber|I_n(t)| &\leq &  C\eps^{-3}(m+1)^2|y_{n+1}(t)-y_n(t)|^{2} \|y_{n-1}(t)\|^{2}
+\eps \|y_n(t)-y_{n-1}(t)\|^{2}
\\ \nonumber  && + \eps \|y_{n-1}(t)\|^2 |y_n(t)-y_{n-1}(t)|^{2} +\eps \|y_{n+1}(t)-y_n(t)\|^{2}.
\end{eqnarray}
The fact that for this subcase $|y_n|_{\xi_t}\leq 3\delta$ and  $|y_{n-1}|_{\xi_t}\leq 3\delta$ show that
\begin{eqnarray}
\nonumber |I_n(t)|
&\leq &  C\eps^{-3}(m+1)^2|y_{n+1}(t)-y_n(t)|^{2} \|y_{n-1}(t)\|^{2}I_{[0,3\delta]}(|y_{n-1}|_{\xi_t})
+\eps \|y_n(t)-y_{n-1}(t)\|^{2}
\\ \label{p-17} && + \eps \|y_{n-1}(t)\|^2 |y_n(t)-y_{n-1}(t)|^{2}I_{[0,3\delta]}(|y_{n-1}|_{\xi_t})   +\eps \|y_{n+1}(t)-y_n(t)\|^{2}.
\end{eqnarray}

\emph{(\uppercase\expandafter{\romannumeral1}-4):   $|y_n(t)|> m+2$ and  $|y_{n-1}(t)|> m+2$.}

For this subcase, note that
\begin{eqnarray}
\label{p-18}
  I_n(t)=0.
\end{eqnarray}

Therefore, (\ref{p-15}), (\ref{p-16}), (\ref{p-17}) and (\ref{p-18}) ensure that for the case (\uppercase\expandafter{\romannumeral1}) the following equality holds, for any $\eps>0$,
\begin{eqnarray}\label{p-19}
 \nonumber I_n(t)
  & \leq&    5\eps   \|y_{n+1}(t)-y_n(t)\|^{2}+2\eps   \|y_n(t)-y_{n-1}(t)\|^{2}\\
&& + \frac{\eps  }{\delta^{3/2}}\|y_n(t)\|^2  |y_n-y_{n-1}|_{\xi_t}^{2} I_{[0,3\delta]}(|y_n|_{\xi_t})
+3\eps |y_n(t)-y_{n-1}(t)|^{2}\Xi_n(t)  \nonumber\\
&&+C\eps^{-3}\big(1+(m+2)^2+\delta^{-1}(m+2)^2\big) |y_{n+1}(t)-y_n(t)|^{2}\Xi_n(t).
\end{eqnarray}
%
%

\textbf{(\uppercase\expandafter{\romannumeral2}):}  \textbf{$|y_n|_{\xi_t}\leq 3\delta$ and  $|y_{n-1}|_{\xi_t}>3\delta$.}

For this case, note that the defintion of $\phi_m$ and $g_\delta$ yields that
\begin{eqnarray*}
 I_n(t)\cdot I_{\{|y_n|_{\xi_t}\leq 3\delta\}}\cdot I_{\{|y_{n-1}|_{\xi_t}> 3\delta\}}
&=&
  I_n(t)\cdot I_{\{|y_n|_{\xi_t}\leq 2\delta\}}\cdot I_{\{|y_{n-1}|_{\xi_t}> 3\delta\}}\cdot I_{\{|y_n(t)|\leq  m+1\}}.
\end{eqnarray*}
For any $t$ such that $I_{\{|y_n|_{\xi_t}\leq 2\delta\}}\cdot I_{\{|y_{n-1}|_{\xi_t}> 3\delta\}}\cdot I_{\{|y_n(t)|\leq  m+1\}}=1$, we have
\begin{eqnarray}\label{p-14}
|y_n-y_{n-1}|_{\xi_t}\geq  \delta,
\end{eqnarray}
and
by (\ref{p-3}), (\ref{p-4}), (\ref{p-5}) and (\ref{eq PB}), we have
\begin{eqnarray}
\nonumber  I_n(t)
  &=&\langle B(y_{n}(t),y_{n+1}(t))\phi_m(|y_n(t)|)g_\delta(|y_n|_{\xi_t}) ,  y_{n+1}(t)-y_n(t) \rangle
\\  \nonumber  &\leq& \big| \langle B(y_{n}(t),y_{n+1}(t)) ,  -y_n(t) \rangle\big|
\\  \nonumber  &=&\big| \langle B(y_{n}(t),y_{n+1}(t)-y_n(t)),  -y_n(t) \rangle\big|
\\  \nonumber  &=&\big| \langle B(y_{n}(t),y_n(t) ),  y_{n+1}(t)-y_n(t) \rangle\big|
\\ \nonumber  &\leq& C |y_n(t)|_Q\cdot \|y_n(t)\|\cdot  |y_{n+1}(t)-y_n(t)|_Q
\\ \nonumber &\leq &  C \|y_{n+1}(t)-y_n(t)\|^{1/2}
\cdot |y_{n+1}(t)-y_n(t)|^{1/2}\cdot \|y_n(t)\|\cdot |y_n(t)|^{1/2}
\|y_n(t)\|^{1/2}.
\end{eqnarray}
 Due to (\ref{p-14}), $|y_n(t)|\leq m+1$, and Young's inequality, we show that for any $\eps,p>0$,
\begin{eqnarray}
\nonumber  I_n(t)
\nonumber   &\leq &  C \|y_{n+1}(t)-y_n(t)\|^{1/2}
\cdot |y_{n+1}(t)-y_n(t)|^{1/2}\cdot \|y_n(t)\|\\  \nonumber
&&\cdot |y_n(t)|^{1/2}
\|y_n(t)\|^{1/2}\cdot |y_n-y_{n-1}|_{\xi_t}\cdot\frac{1}{\delta}
\\  \nonumber  &\leq & C(m+1)^{1/2} \|y_{n+1}(t)-y_n(t)\|^{1/2}\cdot |y_{n+1}(t)-y_n(t)|^{1/2}
\\ \nonumber && \cdot     \|y_n(t)\| \cdot
\|y_n(t)\|^{1/2}\cdot |y_n-y_{n-1}|_{\xi_t}\cdot \eps^{-3/4}\delta^{-p}\cdot \eps^{3/4}\delta^{p-1}
\\  \nonumber  &\leq& C(m+1)^2\eps^{-3} \delta^{-4p}\|y_n(t)\|^2   |y_{n+1}(t)-y_n(t)|^2 \\  \nonumber
&&+\eps \delta^{4(p-1)/3}  \|y_{n+1}(t)-y_n(t)\|^{2/3}\|y_n(t)\|^{4/3} |y_n-y_{n-1}|_{\xi_t}^{4/3}.
\end{eqnarray}
Hence, $|y_n|_{\xi_t}\leq 2\delta<3\delta$ and Young's inequality imply that
\begin{eqnarray}
\nonumber  I_n(t)   &\leq &C(m+1)^2\eps^{-3}  \delta^{-4p} \|y_n(t)\|^2  |y_{n+1}(t)-y_n(t)|^2I_{[0,3\delta]}(|y_{n}|_{\xi_t})
+\eps  \|y_{n+1}(t)-y_n(t)\|^{2}\\
\label{p-20}  && +C\eps \delta^{2(p-1)} \|y_n(t)\|^{2} |y_n-y_{n-1}|_{\xi_t}^{2}I_{[0,3\delta]}(|y_{n}|_{\xi_t}).
\end{eqnarray}

\textbf{(\uppercase\expandafter{\romannumeral3}):}
 \textbf{ $|y_n|_{\xi_t}> 3\delta$ and  $|y_{n-1}|_{\xi_t}\leq 3\delta$.
}

For this case, from the defintion of $\phi_m$ and $g_\delta$, we find
\begin{eqnarray*}
 I_n(t)\cdot I_{\{|y_n|_{\xi_t}> 3\delta\}}\cdot I_{\{|y_{n-1}|_{\xi_t}\leq 3\delta\}}
&=&
  I_n(t)\cdot I_{\{|y_n|_{\xi_t}> 3\delta\}}\cdot I_{\{|y_{n-1}|_{\xi_t}\leq 2\delta\}}\cdot I_{\{|y_{n-1}(t)|\leq  m+1\}}.
\end{eqnarray*}
For any $t$ such that $I_{\{|y_n|_{\xi_t}> 3\delta\}}\cdot I_{\{|y_{n-1}|_{\xi_t}\leq 2\delta\}}\cdot I_{\{|y_{n-1}(t)|\leq  m+1\}}=1$, we have
\begin{eqnarray}\label{p-11}
|y_n-y_{n-1}|_{\xi_t}\geq  \delta.
\end{eqnarray}
This and
(\ref{p-3}), (\ref{p-4}), (\ref{p-5}), $|y_{n-1}(t)|\leq m+1$ and the Young inequality prove that for any $\eps,p>0$,
\begin{eqnarray}
\nonumber  |I_n(t)|
  &=&|\langle B(y_{n-1}(t),y_{n}(t))\phi_m(|y_{n-1}(t)|) g_\delta(|y_{n-1}|_{\xi_t}),  y_{n+1}(t)-y_n(t) \rangle|
\\ \nonumber   &\leq& C |y_{n-1}(t)|_Q \cdot \|y_n(t)\|\cdot |y_{n+1}(t)-y_n(t) |_Q
\\  \nonumber  &\leq & C |y_{n+1}(t)-y_n(t)|^{1/2}\|y_{n+1}(t)-y_n(t)\|^{1/2}\cdot\|y_n(t)\|\cdot
\|y_{n-1}(t)\|^{1/2}|y_{n-1}(t)|^{1/2}
\\ \nonumber  &\leq&  C(m+1)^{1/2}  |y_{n+1}(t)-y_n(t)|^{1/2}\|y_{n+1}(t)-y_n(t)\|^{1/2}\cdot\|y_n(t)\|\cdot
\|y_{n-1}(t)\|^{1/2}  \cdot \delta^{-p}\cdot \delta^p
\\ \nonumber  &\leq& C(m+1)^{2} \delta^{-4p} |y_{n+1}(t)-y_n(t)|^{2}\|y_{n-1}(t)\|^{2}
\\ \nonumber  &&+\delta^{4p/3}\|y_{n+1}(t)-y_n(t)\|^{2/3} \|y_n(t)\|^{4/3}\cdot \eps^{-1/3}\cdot \eps^{1/3}
\\ \nonumber  &\leq &C(m+1)^2  \delta^{-4p}|y_{n+1}(t)-y_n(t)|^{2} \|y_{n-1}(t)\|^{2}+\eps \|y_{n+1}(t)-y_n\|^{2}
+
\eps^{-1/2}\delta^{2p}\|y_n(t)\|^{2}
\\ \nonumber  &\leq& C(m+1)^2   \delta^{-4p}|y_{n+1}(t)-y_n(t)|^{2} \|y_{n-1}(t)\|^{2}+\eps \|y_{n+1}(t)-y_n\|^{2}
\\ \nonumber  &&
+
\eps^{-1/2}\delta^{2p}\|y_n(t)-y_{n-1}(t)\|^{2}
+
\eps^{-1/2}\delta^{2p}\|y_{n-1}(t)\|^{2}
\end{eqnarray}
Hence,
by (\ref{p-11}) and $|y_{n-1}|_{\xi_t}\leq 2\delta<3\delta$, we obtain
\begin{eqnarray}
  \nonumber  &\leq& C(m+1)^2   \delta^{-4p}|y_{n+1}(t)-y_n(t)|^{2} \|y_{n-1}(t)\|^{2}+\eps \|y_{n+1}(t)-y_n\|^{2}
  \\ \nonumber
 && +
\eps^{-1/2}\delta^{2p}\|y_n(t)-y_{n-1}(t)\|^{2}
 +
\eps^{-1/2}\delta^{2p}\|y_{n-1}(t)\|^{2}\cdot \frac{|y_{n}-y_{n-1}|_{\xi_t}^2}{\delta^2}\\
 \nonumber  &=& C(m+1)^2  \delta^{-4p}|y_{n+1}(t)-y_n(t)|^{2} \|y_{n-1}(t)\|^{2}I_{[0,3\delta]}(|y_{n-1}|_{\xi_t})+
\eps^{-1/2}\delta^{2p}\|y_n(t)-y_{n-1}(t)\|^{2}
\\ \label{p-21}  && +\eps \|y_{n+1}(t)-y_n\|^{2}  +
\eps^{-1/2}\delta^{2p-2}\|y_{n-1}(t)\|^{2}\cdot |y_{n}-y_{n-1}|_{\xi_t}^2\cdot I_{[0,3\delta]}(|y_{n-1}|_{\xi_t}).
\end{eqnarray}

\textbf{(\uppercase\expandafter{\romannumeral4}):}
 \textbf{ $|y_n|_{\xi_t}> 3\delta$ and  $|y_{n-1}|_{\xi_t}>3\delta$.
}

In this case, by the defintion of $g_\delta$, we know
\begin{eqnarray}\label{p-22}
  I_n(t)=0.
\end{eqnarray}

\vskip 0.2cm
Summing up the cases (\uppercase\expandafter{\romannumeral1})-
(\uppercase\expandafter{\romannumeral4}),
For any $\eps,p,t>0$,  the following inequality holds,
\begin{eqnarray}\label{Eq In}
 \nonumber I_n(t)
  & \leq&    7\eps   \|y_{n+1}(t)-y_n(t)\|^{2}+\Big(2\eps+\eps^{-1/2}\delta^{2p}\Big)   \|y_n(t)-y_{n-1}(t)\|^{2}
\\ \nonumber && + C\Big(\frac{\eps  }{\delta^{3/2}}+\eps^{-1/2}\delta^{2p-2}+\eps \delta^{2(p-1)}\Big)
|y_n-y_{n-1}|_{\xi_t}^{2} \Xi_n(t)
\\ && +3\eps |y_n(t)-y_{n-1}(t)|^{2}
      \Xi_n(t)
+{C \$_n(t)}|y_{n+1}(t)-y_n(t)|^{2}.
\end{eqnarray}

The proof of Lemma \ref{lemma 2} is complete.
\hfill$\Box$\\

\begin{lemma}\label{lemma 3}
For any $\eps>0$ such that $2-L_2-2\eps>0$ and any {$t_0>0$}, we have
\begin{eqnarray}\label{eq y2-y1 05}
&&\exp\Big(-12\eps\cdot18\delta^2-L_1t_0\Big)\sup_{t\in[0,t_0]}\mathbb{E}\Big(|y_2(t)-y_1(t)|^2\Big)\nonumber\\
&&+
\Big(2-L_2-2\eps\Big)\exp\Big(-12\eps\cdot18\delta^2-L_1t_0\Big)\mathbb{E}\int_0^{t_0}\|y_2(s)-y_1(s)\|^2ds\nonumber\\
&&+
2\exp\Big(-12\eps\cdot18\delta^2-L_1t_0\Big)\mathbb{E}\int_0^{t_0} \Big(6\eps\Xi_2(s)\Big)|y_2(s)-y_1(s)|^2ds\nonumber\\
& \leq&\frac{C}{\eps}m^2\delta^2,
\end{eqnarray}
{where   $C$  is a constant independent of $\eps,\delta, t_0$.}
\end{lemma}

\noindent\emph{Proof of Lemma \ref{lemma 3}.}\quad
Applying $\rm It\hat{o}$'s formula to $\exp\Big(-\int_0^t12\eps\Xi_2(s)+L_1ds\Big)|y_2(t)-y_1(t)|^2$ and
recallint (\ref{eq y1}) and (\ref{XZ-1}), we have
\begin{eqnarray}\label{eq y2-y1 00}
&&\exp\Big(-\int_0^t12\eps\Xi_2(s)+L_1ds\Big)|y_2(t)-y_1(t)|^2
\nonumber\\
&&
+
2\int_0^t \exp\Big(-\int_0^s12\eps\Xi_2(l)+L_1dl\Big)\|y_2(s)-y_1(s)\|^2ds\nonumber\\
&&+
\int_0^t \Big(12\eps\Xi_2(s)+L_1\Big)\exp\Big(-\int_0^s12\eps\Xi_2(l)+L_1dl\Big)|y_2(s)-y_1(s)|^2ds\nonumber\\
&=&
2\int_0^t \exp\Big(-\int_0^s12\eps\Xi_2(l)+L_1dl\Big)
           \langle B(y_{1}(s),y_{2}(s))\phi_m(|y_{1}(s)|) g_\delta(|y_{1}|_{\xi_s}),  y_{2}(s)-y_1(s) \rangle ds\nonumber\\
&&+
2\int_0^t \int_{\mathcal{Z}} \exp\Big(-\int_0^s12\eps\Xi_2(l)+L_1dl\Big)
           \langle G(s,y_2(s-),z)-G(s,y_1(s-),z),
             \nonumber\\
&&\ \ \
             y_{2}(s-)-y_1(s-) \rangle \widetilde{\eta}(dz,ds)\nonumber\\
&&+2\int_0^t \exp\Big(-\int_0^s12\eps\Xi_2(l)+L_1dl\Big)
           \langle \Psi(s,y_2(s))-\Psi(s,y_1(s)),  y_{2}(s)-y_1(s) \rangle dW(s)\nonumber\\
&&+\int_0^t\int_{\mathcal{Z}} \exp\Big(-\int_0^s12\eps\Xi_2(l)+L_1dl\Big)
                              |G(s,y_2(s-),z)-G(s,y_1(s-),z)|^2\eta(dz,ds)\nonumber\\
&&+
\int_0^t \exp\Big(-\int_0^s12\eps\Xi_2(l)+L_1dl\Big)
                |\Psi(s,y_2(s))-\Psi(s,y_1(s))|^2_{\mathcal{L}_2}ds\nonumber\\
&=:&\sum_{i=1}^5 I_i(t).
\end{eqnarray}
According to the fact that $y_1, y_2\in D([0,T],H)\ \mathbb{P}\text{-a.s.}$ and (\ref{eq BJ 05}), Condition \ref{con G}, there exist stopping times $\tau_n\nearrow\infty\ \mathbb{P}\text{-a.s.}$ such that
\begin{eqnarray}\label{eq martingale 1}
\Big\{I_2(t\wedge\tau_n)+I_3(t\wedge\tau_n),\ t\geq0\Big\}\text{ is {an} }\mathbb{F}\text{-martingale}.
\end{eqnarray}
Hence, we obtain
\begin{eqnarray}\label{eq y2-y1 01}
&&\mathbb{E}[\exp\Big(-\int_0^{t\wedge\tau_n}12\eps\Xi_2(s)
+L_1ds\Big)|y_2(t\wedge\tau_n)-y_1(t\wedge\tau_n)|^2]\nonumber\\
&&+
2\mathbb{E}[\int_0^{t\wedge\tau_n} \exp\Big(-\int_0^s12\eps\Xi_2(l)+L_1dl\Big)\|y_2(s)-y_1(s)\|^2ds]\nonumber\\
&&+
\mathbb{E}[\int_0^{t\wedge\tau_n} \Big(12\eps\Xi_2(s)+L_1\Big)\exp\Big(-\int_0^s12\eps\Xi_2(l)+L_1dl\Big)|y_2(s)-y_1(s)|^2ds]\nonumber\\
&=&\mathbb{E}[I_1(t\wedge\tau_n)]+ \mathbb{E}[I_4(t\wedge\tau_n)]+ \mathbb{E}[I_5(t\wedge\tau_n)].
\end{eqnarray}

By Condition \ref{p-2} and (\ref{eq PB}), we have
\begin{eqnarray}\label{eq I1}
&&\mathbb{E}[I_1(t\wedge\tau_n)]\nonumber\\
   &=&
      2\mathbb{E}[\int_0^{t\wedge\tau_n} \exp\Big(-\int_0^s12\eps\Xi_2(l)+L_1dl\Big)
           \langle B(y_{1}(s),y_{1}(s))\phi_m(|y_{1}(s)|) g_\delta(|y_{1}|_{\xi_s}),  \nonumber\\
           &&\ \  y_{2}(s)-y_1(s) \rangle ds]\nonumber\\
  &\leq&
        2\eps \mathbb{E}[\int_0^{t\wedge\tau_n} \exp\Big(-\int_0^s12\eps\Xi_2(l)+L_1dl\Big)\|y_{2}(s)-y_1(s)\|^2ds]\nonumber\\
&&+
        \frac{C}{\eps} \mathbb{E}[\int_0^{t\wedge\tau_n} \exp\Big(-\int_0^s12\eps\Xi_2(l)+L_1dl\Big)
         |y_{1}(s)|^2\|y_{1}(s)\|^2 \phi_m(|y_{1}(s)|) g_\delta(|y_{1}|_{\xi_s})ds]\nonumber\\
  &\leq&
        2\eps \mathbb{E}[\int_0^{t\wedge\tau_n} \exp\Big(-\int_0^s12\eps\Xi_2(l)+L_1dl\Big)\|y_{2}(s)-y_1(s)\|^2ds]
+
        \frac{C}{\eps}m^2\delta^2.
\end{eqnarray}
Condition \ref{con G} shows that
\begin{eqnarray}\label{eq I4I5}
&&\mathbb{E}[I_4(t\wedge\tau_n)]+ \mathbb{E}[I_5(t\wedge\tau_n)]\nonumber\\
   &=&
      \mathbb{E}[\int_0^{t\wedge\tau_n}\int_{\mathcal{Z}} \exp\Big(-\int_0^s12\eps\Xi_2(l)+L_1dl\Big)
                              |G(s,y_2(s),z)-G(s,y_1(s),z)|^2\nu(dz)ds]\nonumber\\
&&+
      \mathbb{E}[\int_0^{t\wedge\tau_n} \exp\Big(-\int_0^s12\eps\Xi_2(l)+L_1dl\Big)
                |\Psi(s,y_2(s))-\Psi(s,y_1(s))|^2_{\mathcal{L}_2}ds]\\
&\leq&
     \mathbb{E}[\int_0^{t\wedge\tau_n} \exp\Big(-\int_0^s12\eps\Xi_2(l)+L_1dl\Big)
                \Big(L_1|y_2(s)-y_1(s)|^2+L_2\|y_2(s)-y_1(s)\|^2\Big)ds].\nonumber
\end{eqnarray}
Combining (\ref{eq y2-y1 01})--(\ref{eq I4I5}), we can yield
\begin{eqnarray}\label{eq y2-y1 02}
&&\mathbb{E}[\exp\Big(-\int_0^{t\wedge\tau_n}12\eps\Xi_2(s)+L_1ds\Big)
|y_2(t\wedge\tau_n)-y_1(t\wedge\tau_n)|^2]\nonumber\\
&&+
\Big(2-L_2-2\eps\Big)\mathbb{E}[\int_0^{t\wedge\tau_n} \exp\Big(-\int_0^s12\eps\Xi_2(l)+L_1dl\Big)\|y_2(s)-y_1(s)\|^2ds]\nonumber\\
&&+
\mathbb{E}[\int_0^{t\wedge\tau_n} \Big(12\eps\Xi_2(s)\Big)\exp\Big(-\int_0^s12\eps\Xi_2(l)+L_1dl\Big)|y_2(s)-y_1(s)|^2ds]\nonumber\\
&\leq&\frac{C}{\eps}m^2\delta^2.
\end{eqnarray}
Taking limit of $n$ tends to infinity assures that
\begin{eqnarray}\label{eq y2-y1 03}
&&\mathbb{E}\Big(\exp\Big(-\int_0^{t}12\eps\Xi_2(s)+L_1ds\Big)|y_2(t)-y_1(t)|^2\Big)\nonumber\\
&&+
\Big(2-L_2-2\eps\Big)\mathbb{E}\int_0^{t} \exp\Big(-\int_0^s12\eps\Xi_2(l)+L_1dl\Big)\|y_2(s)-y_1(s)\|^2ds\nonumber\\
&&+
\mathbb{E}\int_0^{t} \Big(12\eps\Xi_2(s)\Big)\exp\Big(-\int_0^s12\eps\Xi_2(l)+L_1dl\Big)|y_2(s)-y_1(s)|^2ds\nonumber\\
&\leq&\frac{C}{\eps}m^2\delta^2.
\end{eqnarray}
Thanks to the fact that for any $n\geq1$ and $S>0$
\begin{eqnarray*}
\int_0^S \|y_n(t)\|^2I_{[0,3\delta]}(|y_n|_{\xi_t})\dif t\leq 9\delta^2,
\end{eqnarray*}
we obtain
$$
0\leq\int_0^{t\wedge\tau_n}12\eps\Xi_2(s)+L_1ds\leq 12\eps\cdot18\delta^2+L_1t, \forall t\geq 0\text{ and }\forall n.
$$
Applying this inequality to (\ref{eq y2-y1 03}) demonstrates that
\begin{eqnarray}\label{eq y2-y1 04}
&&\exp\Big(-12\eps\cdot18\delta^2-L_1t\Big)\mathbb{E}[|y_2(t)-y_1(t)|^2]\nonumber\\
&&+
\Big(2-L_2-2\eps\Big)\exp\Big(-12\eps\cdot18\delta^2-L_1t\Big)
\mathbb{E}[\int_0^{t}\|y_2(s)-y_1(s)\|^2ds]\nonumber\\
&&+
2\exp\Big(-12\eps\cdot18\delta^2-L_1t\Big)\mathbb{E}[\int_0^{t} \Big(6\eps\Xi_2(s)\Big)|y_2(s)-y_1(s)|^2ds]\nonumber\\
&\leq&\frac{C}{\eps}m^2\delta^2.
\end{eqnarray}

Therefore, we obtain
that for $\eps$ such that $2-L_2-2\eps>0$ and any $t_0>0$,
\begin{eqnarray}\label{eq y2-y1 05}
&&\exp\Big(-12\eps\cdot18\delta^2-L_1t_0\Big)\sup_{t\in[0,t_0]}
\mathbb{E}[|y_2(t)-y_1(t)|^2]\nonumber\\
&&+
\Big(2-L_2-2\eps\Big)\exp\Big(-12\eps\cdot18\delta^2-L_1t_0\Big)
\mathbb{E}[\int_0^{t_0}\|y_2(s)-y_1(s)\|^2ds]\nonumber\\
&&+
2\exp\Big(-12\eps\cdot18\delta^2-L_1t_0\Big)\mathbb{E}[\int_0^{t_0} \Big(6\eps\Xi_2(s)\Big)|y_2(s)-y_1(s)|^2ds]\nonumber\\
&\leq&\frac{C}{\eps}m^2\delta^2.
\end{eqnarray}

The proof of Lemma \ref{lemma 3} is complete.
\hfill$\Box$\\

\vskip 0.2cm

If we can prove Propositions \ref{Th 01} and \ref{thm 02} below, then   (\ref{eq LIp}) follows immediately.

\begin{proposition}\label{Th 01}
{There exist $\delta_0,T_0>0$ independent of  the initial value $u_0$   }such that
\begin{eqnarray}\label{eq yn+1-yn 10-00}
 &&
 \sum_{n=2}^\infty
 \Big(\mathbb{E}[|y_{n+1}(T_0))-y_{n}(T_0)|^2]\Big)^{1/2}+
\sum_{n=2}^\infty
\Big(\mathbb{E}[\int_0^{T_0}\|y_{n+1}(s)-y_{n}(s)\|^2ds]\Big)^{1/2}\nonumber\\
&&+
\sum_{n=2}^\infty
\Big(\mathbb{E}[\int_0^{T_0} \Big(\Xi_{n+1}(s)\Big)|y_{n+1}(s)-y_{n}(s)|^2ds]\Big)^{1/2}
<\infty.
\end{eqnarray}
Here $y_{n+1}$ is the solution of (\ref{XZ-1}) with $\delta$ replaced by $\delta_0$.
\end{proposition}

\noindent\emph{Proof of Proposition \ref{Th 01}.}\quad
Recall the definitions of $\$_n$ and $I_n$ in (\ref{eq Sn}) and (\ref{1-1}), respectively. And set
$$
F_n(s)={\mathcal{C}_0} \$_n(s)+\Xi_{n+1}(s)+L_1,
$$
{where $\mathcal{C}_0$ is a constant will be decided later.}
Applying $\rm It\hat{o}$'s formula to $\exp\Big(-\int_0^t F_n(s)ds\Big)|y_{n+1}(t)-y_n(t)|^2$, we have
\begin{eqnarray}\label{eq yn+1-yn 00}
&&\exp\Big(-\int_0^t F_n(s)ds\Big)|y_{n+1}(t)-y_n(t)|^2\nonumber\\
&&+
2\int_0^t \exp\Big(-\int_0^sF_n(l)dl\Big)\|y_{n+1}(s)-y_n(s)\|^2ds\nonumber\\
&&+
\int_0^t F_n(s)\exp\Big(-\int_0^sF_n(l)dl\Big)|y_{n+1}(s)-y_n(s)|^2ds\nonumber\\
&=&
2\int_0^t \exp\Big(-\int_0^sF_n(l)dl\Big)I_n(s) ds+
2\int_0^t \int_{\mathcal{Z}} \exp\Big(-\int_0^sF_n(l)dl\Big)
\nonumber\\
&&\ \
         \cdot  \langle G(s,y_{n+1}(s-),z)-G(s,y_{n}(s-),z),  y_{n+1}(s-)-y_n(s-) \rangle \widetilde{\eta}(dz,ds)\nonumber\\
&&+2\int_0^t \exp\Big(-\int_0^sF_n(l)dl\Big)
           \langle \Psi(s,y_{n+1}(s))-\Psi(s,y_n(s)),  y_{n+1}(s)-y_n(s) \rangle dW(s)\nonumber\\
&&+\int_0^t\int_{\mathcal{Z}} \exp\Big(-\int_0^sF_n(l)dl\Big)
                              |G(s,y_{n+1}(s-),z)-G(s,y_n(s-),z)|^2\eta(dz,ds)\nonumber\\
&&+
\int_0^t \exp\Big(-\int_0^sF_n(l)dl\Big)
                |\Psi(s,y_{n+1}(s))-\Psi(s,y_n(s))|^2_{\mathcal{L}_2}ds\nonumber\\
&=:&\sum_{i=1}^5 J^n_i(t).
\end{eqnarray}
Since, for any $i\in\mathbb{N}$, $y_i\in D([0,T],H)\ \mathbb{P}\text{-a.s.}$ and by (\ref{eq BJ 05}) and Condition \ref{con G}, for any $N\in\mathbb{N}$, there exist stopping times $\tau^N_\Bbbk\nearrow\infty\ \mathbb{P}\text{-a.s.}$ as $\Bbbk \nearrow\infty$ such that, for any $n\in \{2,3,\cdots, N\}$,
\begin{eqnarray}\label{eq yn+1-yn martingale 1}
\Big\{J^n_2(t\wedge\tau^N_\Bbbk)+J^n_3(t\wedge\tau^N_\Bbbk),\ t\geq0\Big\}\text{ is a }\mathbb{F}\text{-martingale}.
\end{eqnarray}
Hence
\begin{eqnarray}\label{eq yn+1-yn 01}
&&\mathbb{E}[\exp\Big(-\int_0^{t\wedge{\tau^N_\Bbbk}}
F_n(s)ds\Big)|y_{n+1}(t\wedge{\tau^N_\Bbbk})-y_{n}(t\wedge{\tau^N_\Bbbk})|^2]\nonumber\\
&&+
2\mathbb{E}[\int_0^{t\wedge{\tau^N_\Bbbk}} \exp\Big(-\int_0^sF_n(l)dl\Big)\|y_{n+1}(s)-y_{n}(s)\|^2ds]\nonumber\\
&&+
\mathbb{E}[\int_0^{t\wedge{\tau^N_\Bbbk}} \Big(F_n(s)\Big)\exp\Big(-\int_0^sF_n(l)dl\Big)|y_{n+1}(s)-y_{n}(s)|^2ds]\nonumber\\
&=&\mathbb{E}[J^n_1(t\wedge{\tau^N_\Bbbk})]+ \mathbb{E}[J^n_4(t\wedge{\tau^N_\Bbbk})]+ \mathbb{E}[J^n_5(t\wedge{\tau^N_\Bbbk})].
\end{eqnarray}

Similar to (\ref{eq I4I5}), we have
\begin{eqnarray}\label{eq yn+1-yn I4I5}
&&\mathbb{E}[J^n_4(t\wedge{\tau^N_\Bbbk})]+ \mathbb{E}[J^n_5(t\wedge{\tau^N_\Bbbk})]\\
&\leq&
     \mathbb{E}[\int_0^{t\wedge{\tau^N_\Bbbk}} \exp\Big(-\int_0^sF_n(l)dl\Big)
                \Big(L_1|y_{n+1}(s)-y_{n}(s)|^2+L_2\|y_{n+1}(s)-y_{n}(s)\|^2\Big)ds].\nonumber
\end{eqnarray}

Using the fact that
\begin{eqnarray*}
\exp\Big(-\int_0^sF_n(l)dl\Big)\leq 1, \int_0^t\Xi_n(s)ds\leq 18\delta^2\text{ and } |y_n-y_{n-1}|_{\xi_t}^{2}=\int_0^t\|y_n(s)-y_{n-1}(s)\|^{2}ds,
\end{eqnarray*}
it is easy to have the following three estimates:
\begin{eqnarray}\label{eq yn+1-yn J1 1}
\begin{split}
& \mathbb{E}[\int_0^{t\wedge{\tau^N_\Bbbk}}\exp\Big(-\int_0^sF_n(l)dl\Big)
\|y_n(s)-y_{n-1}(s)\|^{2}ds]
\\ & \leq
\mathbb{E}[\int_0^{t\wedge{\tau^N_\Bbbk}}\|y_n(s)-y_{n-1}(s)\|^{2}ds],
\end{split}
\end{eqnarray}
\begin{eqnarray}\label{eq yn+1-yn J1 2}
\begin{split}
&\mathbb{E}[\int_0^{t\wedge{\tau^N_\Bbbk}}\exp\Big(-\int_0^sF_n(l)dl\Big)
|y_n-y_{n-1}|_{\xi_s}^{2} \Xi_n(s)ds]\\
&\leq
\mathbb{E}[|y_n-y_{n-1}|_{\xi_{t\wedge{\tau^N_\Bbbk}}}^{2}\int_0^{t\wedge{\tau^N_\Bbbk}} \Xi_n(s)ds]\\
&\leq
18\delta^2\mathbb{E}[\int_0^{t\wedge{\tau^N_\Bbbk}}\|y_n(s)-y_{n-1}(s)\|^{2}ds],
\end{split}
\end{eqnarray}
\begin{eqnarray}\label{eq yn+1-yn J1 3}
\begin{split}
& \mathbb{E}[\int_0^{t\wedge{\tau^N_\Bbbk}}\exp\Big(-\int_0^sF_n(l)dl\Big)|y_n(s)-y_{n-1}(s)|^{2}
      \Xi_n(s)ds]
      \\
      & \leq
      \mathbb{E}[\int_0^{t\wedge{\tau^N_\Bbbk}}|y_n(s)-y_{n-1}(s)|^{2}
      \Xi_n(s)ds].
      \end{split}
\end{eqnarray}
Throughout our proof, $C$  denotes  the constant  which  appears in (\ref{Eq In}).
Combining (\ref{eq yn+1-yn J1 1})-(\ref{eq yn+1-yn J1 3}) and (\ref{Eq In}) ensure that
\begin{eqnarray}\label{eq yn+1-yn J1}
&&\mathbb{E}[J^n_1(t\wedge{\tau^N_\Bbbk})]\nonumber\\
&\leq&
14\eps\mathbb{E}[\int_0^{t\wedge{\tau^N_\Bbbk}}\exp\Big(-\int_0^sF_n(l)dl\Big)
\|y_{n+1}(s)-y_n(s)\|^{2}ds]\nonumber\\
&&+
\mathbb{C}(\eps,\delta,p)\mathbb{E}[\int_0^{t\wedge{\tau^N_\Bbbk}}\|y_n(s)-y_{n-1}(s)\|^{2}ds]\nonumber\\
&&+
6\eps \mathbb{E}[\int_0^{t\wedge{\tau^N_\Bbbk}}|y_n(s)-y_{n-1}(s)|^{2}
      \Xi_n(s)ds]\nonumber\\
&&+
\mathbb{E}[\int_0^{t\wedge{\tau^N_\Bbbk}}\exp\Big(-\int_0^sF_n(l)dl\Big)
{2C\$_n(s)}|y_{n+1}(s)-y_n(s)|^{2}ds],
\end{eqnarray}
where
\begin{eqnarray}\label{Def C}
\mathbb{C}(\eps,\delta,p)
:=
2\Big(2\eps+\eps^{-1/2}\delta^{2p}+C\Big(\frac{\eps}{\delta^{3/2}}+\eps^{-1/2}\delta^{2p-2}+\eps \delta^{2(p-1)}\Big)\cdot{18\delta^2}\Big).
\end{eqnarray}
Let us set { $\cC_0=2C$.}
Combining (\ref{eq yn+1-yn 01}), (\ref{eq yn+1-yn I4I5}) and (\ref{eq yn+1-yn J1}), we  obtain
 \begin{eqnarray}\label{eq yn+1-yn 02}
&&\mathbb{E}[\exp\Big(-\int_0^{t\wedge{\tau^N_\Bbbk}}F_n(s)ds\Big)|y_{n+1}
(t\wedge{\tau^N_\Bbbk})-y_{n}(t\wedge{\tau^N_\Bbbk})|^2]\nonumber\\
&&+
\Big(2-L_2-14\eps\Big)\mathbb{E}[\int_0^{t\wedge{\tau^N_\Bbbk}} \exp\Big(-\int_0^sF_n(l)dl\Big)\|y_{n+1}(s)-y_{n}(s)\|^2ds]\nonumber\\
&&+
\mathbb{E}[\int_0^{t\wedge{\tau^N_\Bbbk}} \Xi_{n+1}(s)\exp\Big(-\int_0^sF_n(l)dl\Big)|y_{n+1}(s)-y_{n}(s)|^2ds]\nonumber\\
&\leq&
\mathbb{C}(\eps,\delta,p)\mathbb{E}[\int_0^{t\wedge{\tau^N_\Bbbk}}\|y_n(s)-y_{n-1}(s)\|^{2}ds]\nonumber\\
&&+
6\eps \mathbb{E}[\int_0^{t\wedge{\tau^N_\Bbbk}}|y_n(s)-y_{n-1}(s)|^{2}
      \Xi_n(s)ds].
\end{eqnarray}
Using
$
\int_0^t\Xi_n(s)ds\leq 18\delta^2
$
again,
\begin{eqnarray}\label{Def J}
\int_0^sF_n(l)dl
&\leq&
\cC_0 \eps^{-3}\Big(1+(m+2)^2+\delta^{-1}(m+2)^2+(m+1)^2\delta^{-4p}+(m+1)^2 \eps^3 \delta^{-4p}\Big)\cdot18\delta^2
\nonumber\\&&+
18\delta^2
+
L_1s\nonumber\\
&=:&\mathbb{J}(\eps,\delta,p,s),
\end{eqnarray}
which is independent of $n$ and for any $s\in[0,t]$,
$$
\mathbb{J}(\eps,\delta,p,s)\leq \mathbb{J}(\eps,\delta,p,t).
$$
Hence
 \begin{eqnarray}\label{eq yn+1-yn 03}
 &&\exp\Big(-\mathbb{J}(\eps,\delta,p,t)\Big)
 \mathbb{E}[|y_{n+1}(t\wedge{\tau^N_\Bbbk})-y_{n}(t\wedge{\tau^N_\Bbbk})|^2]\nonumber\\
&&+
\Big(2-L_2-14\eps\Big)\exp\Big(-\mathbb{J}(\eps,\delta,p,t)\Big)\mathbb{E}[
\int_0^{t\wedge{\tau^N_\Bbbk}}\|y_{n+1}(s)-y_{n}(s)\|^2ds]\nonumber\\
&&+
\exp\Big(-\mathbb{J}(\eps,\delta,p,t)\Big)\mathbb{E}[\int_0^{t\wedge{\tau^N_\Bbbk}} \Big(\Xi_{n+1}(s)\Big)|y_{n+1}(s)-y_{n}(s)|^2ds]\nonumber\\
 &\leq&
\mathbb{E}[\exp\Big(-\int_0^{t\wedge{\tau^N_\Bbbk}}F_n(s)ds\Big)
|y_{n+1}(t\wedge{\tau^N_\Bbbk})-y_{n}(t\wedge{\tau^N_\Bbbk})|^2]\nonumber\\
&&+
\Big(2-L_2-14\eps\Big)\mathbb{E}[\int_0^{t\wedge{\tau^N_\Bbbk}} \exp\Big(-\int_0^sF_n(l)dl\Big)\|y_{n+1}(s)-y_{n}(s)\|^2ds]\nonumber\\
&&+
\mathbb{E}[\int_0^{t\wedge{\tau^N_\Bbbk}} \Big(\Xi_{n+1}(s)\Big)\exp\Big(-\int_0^sF_n(l)dl\Big)|y_{n+1}(s)-y_{n}(s)|^2ds].
\end{eqnarray}
By (\ref{eq yn+1-yn 02}) and (\ref{eq yn+1-yn 03}), we arrive at
\begin{eqnarray}\label{eq yn+1-yn 03-00}
 &&\exp\Big(-\mathbb{J}(\eps,\delta,p,t)\Big)\mathbb{E}[|y_{n+1}
 (t\wedge{\tau^N_\Bbbk})-y_{n}(t\wedge{\tau^N_\Bbbk})|^2]\nonumber\\
&&+
\Big(2-L_2-14\eps\Big)\exp\Big(-\mathbb{J}(\eps,\delta,p,t)\Big)\mathbb{E}[
\int_0^{t\wedge{\tau^N_\Bbbk}}\|y_{n+1}(s)-y_{n}(s)\|^2ds]\nonumber\\
&&+
\exp\Big(-\mathbb{J}(\eps,\delta,p,t)\Big)\mathbb{E}[\int_0^{t\wedge{\tau^N_\Bbbk}} \Big(\Xi_{n+1}(s)\Big)|y_{n+1}(s)-y_{n}(s)|^2ds]\nonumber\\
&\leq&
\mathbb{C}(\eps,\delta,p)\mathbb{E}[\int_0^{t\wedge{\tau^N_\Bbbk}}\|y_n(s)-y_{n-1}(s)\|^{2}ds]\nonumber\\
&&+
6\eps \mathbb{E}[\int_0^{t\wedge{\tau^N_\Bbbk}}|y_n(s)-y_{n-1}(s)|^{2}
      \Xi_n(s)ds].
\end{eqnarray}
Since for any $a,b,c\geq0$, $a+b+c\geq \frac{1}{3}(\sqrt{a}+\sqrt{b}+\sqrt{c})^{2}$ and $a+b\leq (\sqrt{a}+\sqrt{b})^2$, the above
inequality implies that
\begin{eqnarray}\label{eq yn+1-yn 04}
 &&\frac{\sqrt{3}}{3}\exp\Big(-\mathbb{J}(\eps,\delta,p,t)/2\Big)\Big(\mathbb{E}[
 |y_{n+1}(t\wedge{\tau^N_\Bbbk})-y_{n}(t\wedge{\tau^N_\Bbbk})|^2]\Big)^{1/2}\nonumber\\
&&+
\frac{\sqrt{3}}{3}\Big(2-L_2-14\eps\Big)^{1/2}\exp\Big(-\mathbb{J}(\eps,\delta,p,t)/2\Big)
\Big(\mathbb{E}[\int_0^{t\wedge{\tau^N_\Bbbk}}\|y_{n+1}(s)-y_{n}(s)\|^2ds]\Big)^{1/2}\nonumber\\
&&+
\frac{\sqrt{3}}{3}\exp\Big(-\mathbb{J}(\eps,\delta,p,t)/2\Big)
\Big(\mathbb{E}[\int_0^{t\wedge{\tau^N_\Bbbk}} \Big(\Xi_{n+1}(s)\Big)|y_{n+1}(s)-y_{n}(s)|^2ds]\Big)^{1/2}\nonumber\\
&\leq&
\mathbb{C}(\eps,\delta,p)^{1/2}\Big(\mathbb{E}
[\int_0^{t\wedge{\tau^N_\Bbbk}}\|y_n(s)-y_{n-1}(s)\|^{2}ds]\Big)^{1/2}\nonumber\\
&&+
(6\eps)^{1/2} \Big(\mathbb{E}[\int_0^{t\wedge{\tau^N_\Bbbk}}|y_n(s)-y_{n-1}(s)|^{2}
      \Xi_n(s)ds]\Big)^{1/2}.
\end{eqnarray}
Summing $n$ from 2 to $N$, we obtain
\begin{eqnarray}\label{eq yn+1-yn 05}
 &&\frac{\sqrt{3}}{3}\exp\Big(-\mathbb{J}(\eps,\delta,p,t)/2\Big)
 \sum_{n=2}^N
 \Big(\mathbb{E}[|y_{n+1}(t\wedge{\tau^N_\Bbbk})-y_{n}(t\wedge{\tau^N_\Bbbk})|^2]\Big)^{1/2}\nonumber\\
&&+
\frac{\sqrt{3}}{3}\Big(2-L_2-14\eps\Big)^{1/2}\exp\Big(-\mathbb{J}(\eps,\delta,p,t)/2\Big)
\sum_{n=2}^N
\Big(\mathbb{E}[\int_0^{t\wedge{\tau^N_\Bbbk}}\|y_{n+1}(s)-y_{n}(s)\|^2ds]\Big)^{1/2}\nonumber\\
&&+
\frac{\sqrt{3}}{3}\exp\Big(-\mathbb{J}(\eps,\delta,p,t)/2\Big)
\sum_{n=2}^N
\Big(\mathbb{E}[\int_0^{t\wedge{\tau^N_\Bbbk}} \Big(\Xi_{n+1}(s)\Big)|y_{n+1}(s)-y_{n}(s)|^2ds]\Big)^{1/2}\nonumber\\
&\leq&
\mathbb{C}(\eps,\delta,p)^{1/2}\sum_{n=2}^N\Big(\mathbb{E}[\int_0^{t\wedge{\tau^N_\Bbbk}}
\|y_n(s)-y_{n-1}(s)\|^{2}ds]\Big)^{1/2}\nonumber\\
&&+
(6\eps)^{1/2} \sum_{n=2}^N\Big(\mathbb{E}[\int_0^{t\wedge{\tau^N_\Bbbk}}|y_n(s)-y_{n-1}(s)|^{2}
      \Xi_n(s)ds]\Big)^{1/2}\nonumber\\
&=&
\mathbb{C}(\eps,\delta,p)^{1/2}\sum_{n=1}^{N-1}\Big(\mathbb{E}[\int_0^{t\wedge{\tau^N_\Bbbk}}
\|y_{n+1}(s)-y_{n}(s)\|^{2}ds]\Big)^{1/2}\nonumber\\
&&+
(6\eps)^{1/2} \sum_{n=1}^{N-1}\Big(\mathbb{E}[\int_0^{t\wedge{\tau^N_\Bbbk}}|y_{n+1}(s)-y_{n}(s)|^{2}
      \Xi_{n+1}(s)ds]\Big)^{1/2}.
\end{eqnarray}
Arranging the inequality above gives that
\begin{eqnarray}\label{eq yn+1-yn 06}
 &&\frac{\sqrt{3}}{3}\exp\Big(-\mathbb{J}(\eps,\delta,p,t)/2\Big)
 \sum_{n=2}^N
 \Big(\mathbb{E}[|y_{n+1}(t\wedge{\tau^N_\Bbbk})-y_{n}(t\wedge{\tau^N_\Bbbk})|^2]\Big)^{1/2}\nonumber\\
&&+
\Big[\frac{\sqrt{3}}{3}\Big(2-L_2-14\eps\Big)^{1/2}\exp\Big(-\mathbb{J}(\eps,\delta,p,t)/2\Big)-\mathbb{C}(\eps,\delta,p)^{1/2}\Big]
\nonumber\\
&&\ \  \cdot
\sum_{n=2}^N
\Big(\mathbb{E}[\int_0^{t\wedge{\tau^N_\Bbbk}}\|y_{n+1}(s)-y_{n}(s)\|^2ds]\Big)^{1/2}\nonumber\\
&&+
\Big[\frac{\sqrt{3}}{3}\exp\Big(-\mathbb{J}(\eps,\delta,p,t)/2\Big)-(6\eps)^{1/2}\Big]
\nonumber\\
&&\ \  \cdot
\sum_{n=2}^N
\Big(\mathbb{E}[\int_0^{t\wedge{\tau^N_\Bbbk}} \Big(\Xi_{n+1}(s)\Big)|y_{n+1}(s)-y_{n}(s)|^2ds]\Big)^{1/2}\nonumber\\
&\leq&
\mathbb{C}(\eps,\delta,p)^{1/2}\Big(\mathbb{E}[\int_0^{t\wedge{\tau^N_\Bbbk}}\|y_{2}(s)-y_{1}(s)\|^{2}
ds]\Big)^{1/2}
\nonumber\\
&&
+
(6\eps)^{1/2}\Big(\mathbb{E}[\int_0^{t\wedge{\tau^N_\Bbbk}}|y_{2}(s)-y_{1}(s)|^{2}
      \Xi_{2}(s)ds]\Big)^{1/2}.
\end{eqnarray}
{Let $p=1/4$ and $\eps=\delta^{\frac{1}{4}}$.  Noting  that $L_2<2$(see Condition \ref{con G}), by the definitions of $\mathbb{C}$ and $\mathbb{J}$ in (\ref{Def C}) and (\ref{Def J}), there exist positive constants $\delta_0,\ T_0,\ \eta_i,\ i=0,1,2\cdots,5$, and denote $\eps_0=\delta_0^{\frac{1}{4}}$, such that}
\begin{eqnarray}\label{eq yn+1-yn 07}
&& \frac{\sqrt{3}}{3}\exp\Big(-\mathbb{J}(\eps_0,\delta_0,p,T_0)/2\Big)\geq \eta_0,\\
&& \Big[\frac{\sqrt{3}}{3}\Big(2-L_2-14\eps_0\Big)^{1/2}\exp\Big(-\mathbb{J}(\eps_0,\delta_0,p,T_0)/2\Big)-\mathbb{C}(\eps_0,\delta_0,p)^{1/2}\Big]
\geq \eta_1,\\
&& \Big[\frac{\sqrt{3}}{3}\exp\Big(-\mathbb{J}(\eps_0,\delta_0,p,T_0)/2\Big)-(6\eps_0)^{1/2}\Big]
\geq \eta_2,\\
&& \Big(2-L_2-2\eps_0\Big)\exp\Big(-12\eps_0\cdot18\delta_0^2-L_1T_0\Big)
\geq \eta_3,\\
&&2\exp\Big(-12\eps_0\cdot18\delta_0^2-L_1T_0\Big)
\geq \eta_4.
\end{eqnarray}
Then, by (\ref{eq y2-y1 05}) and (\ref{eq yn+1-yn 06}), we yield
\begin{eqnarray}\label{eq yn+1-yn 08}
&&\eta_3\mathbb{E}[\int_0^{T_0}\|y_2(s)-y_1(s)\|^2ds]
+
\eta_4\mathbb{E}[\int_0^{T_0} \Big(6\eps_0\Xi_2(s)\Big)|y_2(s)-y_1(s)|^2ds]\nonumber\\
&\leq&\frac{C}{\eps_0}m^2\delta_0^2,
\end{eqnarray}
and
\begin{eqnarray}\label{eq yn+1-yn 09}
 &&\eta_0
 \sum_{n=2}^N
 \Big(\mathbb{E}[|y_{n+1}(T_0)\wedge{\tau^N_\Bbbk})-y_{n}(T_0\wedge{\tau^N_\Bbbk})|^2]\Big)^{1/2}\nonumber\\
&&+
\eta_1
\sum_{n=2}^N
\Big(\mathbb{E}[\int_0^{T_0\wedge{\tau^N_\Bbbk}}\|y_{n+1}(s)-y_{n}(s)\|^2ds]\Big)^{1/2}\nonumber\\
&&+
\eta_2
\sum_{n=2}^N
\Big(\mathbb{E}[\int_0^{T_0\wedge{\tau^N_\Bbbk}} \Big(\Xi_{n+1}(s)\Big)|y_{n+1}(s)-y_{n}(s)|^2ds]\Big)^{1/2}\nonumber\\
&\leq&
\mathbb{C}(\eps_0,\delta_0,p)^{1/2}\Big(\mathbb{E}[\int_0^{T_0\wedge{\tau^N_\Bbbk}}
\|y_{2}(s)-y_{1}(s)\|^{2}ds]\Big)^{1/2}\nonumber\\
&&+
(6\eps_0)^{1/2}\Big(\mathbb{E}[\int_0^{T_0\wedge{\tau^N_\Bbbk}}|y_{2}(s)-y_{1}(s)|^{2}
      \Xi_{2}(s)ds]\Big)^{1/2}.
\end{eqnarray}
Notice that here $y_{n},\ n\geq 1$ now is the solution of (\ref{XZ-1}) with $\delta$ replaced by $\delta_0$.

Let $\Bbbk\nearrow\infty$ firstly and then let $N\nearrow\infty$ in (\ref{eq yn+1-yn 09}), using (\ref{eq yn+1-yn 08}) to obtain
\begin{eqnarray}\label{eq yn+1-yn 10}
 &&\eta_0
 \sum_{n=2}^\infty
 \Big(\mathbb{E}[ |y_{n+1}(T_0))-y_{n}(T_0)|^2] \Big)^{1/2}\nonumber\\
&&+
\eta_1
\sum_{n=2}^\infty
\Big(\mathbb{E}[\int_0^{T_0}\|y_{n+1}(s)-y_{n}(s)\|^2ds]\Big)^{1/2}\nonumber\\
&&+
\eta_2
\sum_{n=2}^\infty
\Big(\mathbb{E}[\int_0^{T_0} \Big(\Xi_{n+1}(s)\Big)|y_{n+1}(s)-y_{n}(s)|^2ds]\Big)^{1/2}\nonumber\\
&\leq&
\mathbb{C}(\eps_0,\delta_0,p)^{1/2}\Big[\mathbb{E}(\int_0^{T_0}\|y_{2}(s)-y_{1}(s)\|^{2}ds]\Big)^{1/2}\nonumber\\
&&+
(6\eps_0)^{1/2}\Big(\mathbb{E}[\int_0^{T_0}|y_{2}(s)-y_{1}(s)|^{2}
      \Xi_{2}(s)ds]\Big)^{1/2}\nonumber\\
&<&\infty.
\end{eqnarray}

The proof of Proposition \ref{Th 01} is complete.
\hfill$\Box$\\

\vskip 0.2cm

\begin{proposition}\label{thm 02}
Assume that $y_{n}, n\geq 1$ is the solution of (\ref{XZ-1}) with $\delta$ replaced by $\delta_0$ and {$\delta_0,T_0>0$  are the constants appearing in Proposition \ref{Th 01}. We have
$$\sum_{n=2}^\infty\mathbb{E}[\sup_{t\in[0,T_0]}|y_{n+1}(t))-y_{n}(t)|]<\infty$$.}
\end{proposition}
\noindent\emph{Proof of Proposition \ref{thm 02}.}\quad
Now, for any $\lambda>0$, define a stopping time
$$
\tau^n_\lambda:=\inf\{s\geq0:\ |y_{n+1}(s)-y_n(s)|\geq \lambda\}.
$$
Since $\{|y_{n+1}(s)-y_n(s)|,\ s\geq0\}$ is $\rm c\grave{a}dl\grave{a}g$, we have $|y_{n+1}(\tau^n_\lambda)-y_n(\tau^n_\lambda)|\geq \lambda$
and
\begin{eqnarray*}
\lambda^2\mathbb{P}\Big(\sup_{t\in[0,T_0]}|y_{n+1}(t)-y_n(t)|>\lambda\Big)
&\leq&
\lambda^2\mathbb{P}\Big(\tau^n_\lambda\leq T_0\Big)\\
&\leq&
\mathbb{E}[|y_{n+1}(T_0\wedge\tau^n_\lambda))-y_{n}
(T_0\wedge\tau^n_\lambda)|^2]  \\
&=:&\kappa_n,
\end{eqnarray*}
and
\begin{eqnarray}\label{eq star 0}
\mathbb{E}[ \sup_{t\in[0,T_0]}|y_{n+1}(t))-y_{n}(t)|]
&=&
\int_0^\infty \mathbb{P}\Big(\sup_{t\in[0,T_0]}|y_{n+1}(t)-y_n(t)|>\lambda\Big)d\lambda\nonumber\\
&\leq&
\int_0^\infty (\lambda^{-2}\kappa_n)\wedge 1d\lambda\nonumber\\
&=&
2\kappa_n^{1/2}.
\end{eqnarray}

Now we estimate $\sum_{n=2}^\infty \kappa_n^{1/2}$.

Using suitable stopping time technique, similar to (\ref{eq yn+1-yn 04}), we can obtain
\begin{eqnarray*}
 &&\frac{\sqrt{3}}{3}\exp\Big(-\mathbb{J}(\eps_0,\delta_0,p,T_0)/2\Big)
 \Big(\mathbb{E}[|y_{n+1}(T_0\wedge\tau^n_\lambda\wedge{\tau^N_\Bbbk})
 -y_{n}({T_0}\wedge\tau^n_\lambda\wedge{\tau^N_\Bbbk})|^2]\Big)^{1/2}\nonumber\\
&\leq&
\mathbb{C}(\eps_0,\delta_0,p)^{1/2}\Big(\mathbb{E}[\int_0^{T_0\wedge\tau^n_\lambda\wedge{\tau^N_\Bbbk}}
\|y_n(s)-y_{n-1}(s)\|^{2}ds]\Big)^{1/2}\nonumber\\
&&+
(6\eps_0)^{1/2} \Big(\mathbb{E}[\int_0^{T_0\wedge\tau^n_\lambda\wedge{\tau^N_\Bbbk}}|y_n(s)-y_{n-1}(s)|^{2}
      \Xi_n(s)ds\Big)\Big]^{1/2}\nonumber\\
&\leq&
\mathbb{C}(\eps_0,\delta_0,p)^{1/2}\Big(\mathbb{E}[\int_0^{T_0}\|y_n(s)-y_{n-1}(s)\|^{2}ds]\Big)^{1/2}\nonumber\\
&&+
(6\eps_0)^{1/2} \Big(\mathbb{E}[\int_0^{T_0}|y_n(s)-y_{n-1}(s)|^{2}
      \Xi_n(s)ds]\Big)^{1/2}.
\end{eqnarray*}
In the above inequality, let $\Bbbk\nearrow\infty$ firstly, and then sum $n$ from 2 to $\infty$, we get
\begin{eqnarray}\label{eq star}
 &&\frac{\sqrt{3}}{3}\exp\Big(-\mathbb{J}(\eps_0,\delta_0,p,T_0)/2\Big)\sum_{n=2}^\infty \kappa_n^{1/2}
 \nonumber\\
&\leq&
\mathbb{C}(\eps_0,\delta_0,p)^{1/2}\sum_{n=2}^\infty \Big(\mathbb{E}[\int_0^{T_0}\|y_n(s)-y_{n-1}(s)\|^{2}ds]\Big)^{1/2}\nonumber\\
&&+
(6\eps_0)^{1/2} \sum_{n=2}^\infty \Big(\mathbb{E}[\int_0^{T_0}|y_n(s)-y_{n-1}(s)|^{2}
      \Xi_n(s)ds]\Big)^{1/2}\nonumber\\
&<&
\infty,
\end{eqnarray}
(\ref{eq yn+1-yn 10}) has been used to get the last inequality.

Combining (\ref{eq star 0}) and (\ref{eq star}) yield that
\begin{eqnarray}\label{eq Lip C}
\sum_{n=2}^\infty\mathbb{E}[\sup_{t\in[0,T_0]}|y_{n+1}(t)-y_{n}(t)|]<\infty.
\end{eqnarray}

The proof of Proposition \ref{thm 02} is complete.
\hfill$\Box$\\
\vskip 0.2cm

From Propositions \ref{Th 01} and \ref{thm 02}, we now prove the following result, which implies the local existence of (\ref{p-1}).

\begin{proposition}\label{thm 03}
For any $T>0$ and $m>0$, there exists a solution to the following  equation on $[0,T].$
\begin{eqnarray}\label{p-27}
\left\{
\begin{split}
& d u(t)  +\mathcal{A}u(t)dt
+B(u(t),u(t))\phi_m(|u(t)|) dt
\\ & = f(t)dt + \int_{\mathcal{Z}}G(t,u(t-),z)\widetilde{\eta}(dz,dt)+\Psi(t,u(t))d W(t),
\\
& u(0)=u_0.
\end{split}
\right.
\end{eqnarray}
\end{proposition}

\noindent\emph{Proof of Proposition \ref{thm 03}.}\quad
Propositions \ref{Th 01} and \ref{thm 02} imply that, for any fixed $m\in \mN$, there exists $T_0>0$,   $\delta_0>0$, $Y_1\in \Upsilon_{T_0}\ \mathbb{P}$-a.s. and a subsequence of $\{y_n(t),t\in[0,T_0]\}_{n\in\mathbb{N}}$, denoted by $\{y_{n_k}(t),t\in[0,T_0]\}_{k\in\mathbb{N}}$, such that,
\begin{eqnarray}
\lim_{k\nearrow\infty}\sup_{t\in[0,T_0]}|Y_1(t))-y_{n_k}(t)|=0,\text{ and }\lim_{k\nearrow\infty}\int_0^{T_0}\|Y_1(t))-y_{n_k}(t)\|^2dt=0,\ \mathbb{P}\text{-a.s..}
\end{eqnarray}
And then it is not difficulty  to prove that $\{Y_1(t),t\in[0,T_0]\}$ is a solution of the following SPDE:
%
%
%
%
%
%
%
%
\begin{eqnarray}\label{eq 03}
\left\{
\begin{split}
& d y +\mathcal{A}ydt
+B(y(t),y(t))\phi_m(|y(t)|)
g_{\delta_0}(|y|_{\xi_t})dt
\\ & = f(t)dt + \int_{\mathcal{Z}}G(t,y(t-),z)\widetilde{\eta}(dz,dt)+\Psi(t,y(t))d W(t),
\\
& y(0)=u_0.
\end{split}
\right.
\end{eqnarray}
Let $\tau_0=0$, and
\begin{eqnarray*}
  \tau_1:=\inf\{t\geq 0, |Y_1|_{\xi_t}\geq \delta_0\}\wedge T_0.
\end{eqnarray*}

By induction, consider the following time-inhomogeneous equation, for $i\geq 1$,
\begin{eqnarray}\label{eq 04}
\left\{
\begin{split}
& y(t)=Y_{i}(t),\ t\in[0,\tau_i],\\
& y(t)+\int_{\tau_{i}}^{t}\mathcal{A}y(s)ds
+\int_{\tau_{i}}^{t} B(y(s),y(s))\phi_m(|y(s)|)
g_{\delta_0}(|y|_{\xi^i_s})ds
\\ &= Y_i(\tau_{i})+\int_{\tau_{i}}^{t} f(s)ds + \int_{\tau_{i}}^{t} \int_{\mathcal{Z}}G(s,y(s-),z)\widetilde{\eta}(dz,ds)+\int_{\tau_{i}}^{t} \Psi(s,y(s))d W(s),  ~t\in[\tau_{i},\tau_i+T_0],\\
& \tau_{i+1}:=\inf\{t\geq \tau_i:\ \int_{\tau_i}^t\|Y_{i+1}(s)\|^2ds\geq\delta\}\wedge (\tau_i+T_0).
\end{split}
\right.
\end{eqnarray}
Here for any $h(\omega)\in L^2_{loc}([0,\infty),V)$, $|h(\omega)|_{\xi^i_s}:=\Big(\int_{\tau_i(\omega)}^s\|h(l,\omega)\|^2dl\Big)^{1/2},\ \forall s\geq\tau_i(\omega)$.

Note that
 the value of  $T_0$ is independent of the initial value.
Using similar arguments as proving $\{Y_1(t),t\in[0,T_0]\}$ is a solution of (\ref{eq 03}), there exists $\{Y_{i+1}(t), t\in [0,\tau_i+T_0]\}$ which is a solution of (\ref{eq 04}). It is not difficulty to see that $\{Y_n,\ {\tau_n}\}_{n\in\mathbb{N}}$
satisfying
\begin{itemize}
\item $Y_n(\omega)\in \Upsilon_{{\tau_n}(\omega)}$ $\mathbb{P}$-a.s.,

\item $Y_n(t\wedge{\tau_n})\in\mathcal{F}_t$, $\forall t\geq 0$,


\item $0=\tau_0\leq\tau_1\leq\tau_2\cdots\leq {\tau_n}\leq\tau_{n+1}\leq \cdots$,

\item $Y_{n+1}(t)=Y_n(t)$, $t\in[0,{\tau_n}]$  $\mathbb{P}$-a.s.,

\item $Y_{n}$ is a solution of (\ref{p-27}) on $[0,\tau_{n}]$.
\end{itemize}

Define
\begin{eqnarray*}
  u(t)=Y_n(t),\ \ \ t\in[0, {\tau_n}].
\end{eqnarray*}
Let
\begin{eqnarray*}
  \tau_{\max}=\lim_{n\rightarrow \infty}{\tau_n}.
\end{eqnarray*}
Then $u$ is a solution to the equation (\ref{p-27}) on $[0,\tau_{\max})$, and
\begin{eqnarray}\label{eq 05}
\int_0^{\tau_{\max}}\|u(s)\|^2ds=\infty,\ \ \ \ \text{ on } \{\omega,\ \tau_{\max}<\infty\}\ \mathbb{P}\text{-a.s.}
\end{eqnarray}
 We explain the equality above  in detail.
By the definition of  $\tau_{i+1},$  for each $i$,
either
$$\int_{\tau_i}^{\tau_{i+1}}\|Y_{i+1}(s)\|^2ds=\delta$$ or $\tau_{i+1}=\tau_i+T_0$ holds.
Since $\tau_{\max}(\omega)<\infty,$  there  are only  finitely many $i$ such that $\tau_{i+1}=\tau_i+T_0$, and
there are
infinite many $i$ with $\int_{\tau_i}^{\tau_{i+1}}\|Y_{i+1}(s)\|^2ds=\delta$,
which  implies (\ref{eq 05}).

It remains to show that
$\mP(\tau_{\max}\geq T)=1$, $\forall T\geq 0$.

By $\rm It\hat{o}$'s formula to $|u(t)|^2$, we have
\begin{eqnarray*}
&&|u(t)|^2+2\int_0^t\|u(s)\|^2ds\\
&=&
   |u_0|^2
   +
   2\int_0^t\langle f(s),u(s)\rangle ds  +
   2\int_0^t\int_{\mathcal{Z}}\langle G(s,u(s-),z),u(s-)\rangle \widetilde{\eta}(dz,ds)\\
   &&+
    \int_0^t\int_{\mathcal{Z}}|G(s,u(s-),z)|^2 \eta(dzds)\\
    &&
    +2\int_0^t\langle u(s), \Psi(s,u(s))d W(s)\rangle+\int_0^t\|\Psi(s,u(s))\|_{\cL_2}^2d s,\  ~~ \forall t\in[0,T\wedge{\tau_n}].
\end{eqnarray*}
Using suitable stopping time technique and Condition \ref{con G},
we indeed have
\begin{eqnarray}
&&   \mathbb{E}[|u(T\wedge{\tau_n})|^2]
       +
     2\mathbb{E}[\int_0^{T\wedge{\tau_n}}\|u(s)\|^2ds]\nonumber\\
&\leq&
    \mE[  |u_0|^2]
     +
     2\mE [\int_0^{T\wedge{\tau_n}}\langle f(s),u(s)\rangle ds]
       +
       \mathbb{E}[\int_0^{T\wedge{\tau_n}}\int_{\mathcal{Z}} |G(s,u(s),z)|^2 \nu(dz)ds]\nonumber\\
       &&+
       \mE[\int_0^{T\wedge {\tau_n}}\|\Psi(s,u(s))\|_{\cL_2}^2d s]
       \nonumber\\
&\leq&
  \mE [ |u_0|^2]
   +
     \frac{1}{\eps} \int_0^T\|f(s)\|^2_{V'}ds+\eps\mathbb{E}[\int_0^{T\wedge{\tau_n}}\|u(s)\|^2ds]\nonumber\\
     &&+
      L_4\mathbb{E}[\int_0^{T\wedge {\tau_n}}|u(t)|^2dt]
   +
      L_5\mathbb{E}[\int_0^{T\wedge {\tau_n}}\|u(t)\|^2dt]
   +
   L_3T.\nonumber
\end{eqnarray}
Hence
\begin{eqnarray}
&&   \mathbb{E}[|u(T\wedge{\tau_n})|^2]
       +
     (2-L_5-\eps)\mathbb{E}[\int_0^{T\wedge{\tau_n}}\|u(s)\|^2ds]\nonumber\\
&\leq&
  \mE [ |u_0|^2]
   +
     \frac{1}{\eps} \int_0^T\|f(s)\|^2_{V'}ds
     +
      L_4\Big(\int_0^{T}\mathbb{E}[|u(t\wedge {\tau_n})|^2] dt\Big)
   +
   L_3T.\nonumber
\end{eqnarray}

Let $\eps=\frac{2-L_5}{2}$, and by  Gronwall's lemma,
\begin{eqnarray*}
\mathbb{E}[|u(T\wedge{\tau_n}) |^2]
       +
     \mathbb{E}[\int_0^{T\wedge{\tau_n}} \|u(s)\|^2ds]
     \leq
     C_T(1+\mathbb{E}[|u_0|^2 ])<\infty.
\end{eqnarray*}
Taking $n\uparrow\infty$, we have
\begin{eqnarray}\label{eq TT}
\mathbb{E}[ |u(T\wedge{\tau_{\max}}) |^2]
       +
     \mathbb{E}[ \int_0^{T\wedge\tau_{\max}} \|u(s)\|^2ds]
     \leq
     C_T(1+\mathbb{E}[|u_0|^2] )<\infty.
\end{eqnarray}
The above inequality and (\ref{eq 05}) imply that
\begin{eqnarray}\label{eq 06}
\mathbb{P}(\tau_{\max}\geq T)=1.
\end{eqnarray}

The proof of Proposition \ref{thm 03} is complete.
\hfill$\Box$\\

Now we prove Theorem \ref{thm solution}.

\noindent\emph{Proof of Theorem \ref{thm solution}.}\quad
For the uniqueness, we refer to \cite{liuwei} or \cite{BHZ}. In the following, we will prove the global existence.

By the results in Proposition \ref{thm 03}, for any $m\in \mN$,  let $U_m$ be a solution to the following equation

\begin{eqnarray}\label{p-29}
\left\{
\begin{split}
& d u(t)  +\mathcal{A}u(t)dt
+B(u(t),u(t))\phi_m(|u(t)|) dt
\\ & = f(t)dt + \int_{\mathcal{Z}}G(t,u(t-),z)\widetilde{\eta}(dz,dt)+\Psi(t,u(t))d W(t),
\\
& u(0)=u_0.
\end{split}
\right.
\end{eqnarray}

Define
\begin{eqnarray*}
  \sigma_m=\inf\{t\geq 0, |U_m(t)|\geq m\}.
\end{eqnarray*}
Then $\{U_m(t),\ t\in[0,\sigma_m]\}$ is a local solution of (\ref{p-1}).
By the uniqueness,
\begin{eqnarray*}
  U_{m+1}(t)=U_m(t), \text{ on } t\leq \sigma_m,
\end{eqnarray*}
and
\begin{eqnarray*}
  \sigma_m\leq \sigma_{m+1}.
\end{eqnarray*}

Let
\begin{eqnarray*}
  u(t)=U_m(t), \text{ if } t\leq \sigma_m,\text{  and   }\sigma_{\max}=\lim_{ m \rightarrow \infty}\sigma_m.
\end{eqnarray*}

It is obviously that  $u$ is a solution to  (\ref{p-1}) on $[0,\sigma_{\max})$, and
$$
\lim_{m\nearrow\infty}|u(\sigma_m)|=\infty,\ \ \text{on }\{\omega,\sigma_{\max}<\infty\}\ \ \mP\text{-a.s.}.
$$

Using the similar arguments as proving (\ref{eq TT}), we can get, for any $T>0$,
\begin{eqnarray*}
\mathbb{E}[  |u(T\wedge{\sigma_{\max}}) |^2 ]
       +
     \mathbb{E}[  \int_0^{T\wedge\sigma_{\max}} \|u(s)\|^2ds ]
     \leq
     C_T(1+|u_0|^2 ).
\end{eqnarray*}
which implies that
$$
\mathbb{P}(\sigma_{max}\leq T)=0.
$$
and hence
$$
\mathbb{P}(\sigma_{max}=\infty)=1.
$$

The proof of Theorem \ref{thm solution} is complete.
\hfill$\Box$\\

\noindent$\bf{Acknowledgement}$

The authors are very grateful to Professors
Ping Cao and Yong Liu   for their valuable suggestions.

\end{document}